%% file: PsiDiagrams.tex
\documentclass [12pt]{amsart}
\usepackage {graphicx,color}
\usepackage {amsmath,amssymb,latexsym}
\usepackage {bm}
\usepackage {curves}
\usepackage [dvips,centering,includehead,width=15.1cm,height=22.7cm]{geometry}

\input {definitions.tex}

\title {Psi-floor diagrams and a Caporaso-Harris type recursion}
\author {Florian Block \and Andreas Gathmann \and Hannah Markwig}
\address {Florian Block, Department of Mathematics, University of
  Michigan, Ann Arbor, MI 48109, USA} \email{blockf@umich.edu}
\address {Andreas Gathmann, Fachbereich Mathematik, Technische Universit\"at
  Kaiserslautern, Postfach 3049, 67653 Kaiserslautern, Germany}
  \email{andreas@mathematik.uni-kl.de}
\address {Hannah Markwig, Universit\"at des Saarlandes, Fachrichtung
  Mathematik, Postfach 151150, 66041 Saarbr\"ucken, Germany}
  \email{hannah@uni-math.gwdg.de}
\thanks {\emph {2010 Mathematics Subject Classification:} 14T05, 14N35, 51M20}
\keywords {Tropical geometry, enumerative geometry, Gromov-Witten theory}

\begin {document}

\begin {abstract}
  Floor diagrams are combinatorial objects which organize the count of tropical
  plane curves satisfying point conditions. In this paper we introduce
  Psi-floor diagrams which count tropical curves satisfying not only point
  conditions but also conditions given by Psi-classes (together with points).
  We then generalize our definition to relative Psi-floor diagrams and prove a
  Caporaso-Harris type formula for the corresponding numbers. This formula is
  shown to coincide with the classical Caporaso-Harris formula for relative
  plane descendant Gromov-Witten invariants. As a consequence, we can conclude
  that in our case relative descendant Gromov-Witten invariants equal their
  tropical counterparts.
\end {abstract}

\maketitle

\input {introduction.tex} 
\input {classnumbers.tex} 
\input {tropnumbers.tex} 
\input {psifloor.tex} 

\input {biblio.tex}

\end {document}

%% file: definitions.tex
\numberwithin{equation}{section}

\newtheorem{theorem}{Theorem}[section]
\newtheorem{lemma}[theorem]{Lemma}

\newtheorem{notation}[theorem]{Notation}

\newtheorem{corollary}[theorem]{Corollary}
\newtheorem{proposition}[theorem]{Proposition}

\theoremstyle{definition}
\newtheorem{definition}[theorem]{Definition}
\newtheorem{example}[theorem]{Example}
\newtheorem{remark}[theorem]{Remark}
\newtheorem {construction}[theorem]{Construction}

\newcommand{\C}{{\mathcal C}}

\newcommand{\D}{{\mathcal D}}

 \newcommand{\RR}{{\mathbb R}}
\newcommand{\ZZ}{{\mathbb Z}} \newcommand{\NN}{{\mathbb N}}
\newcommand{\QQ}{{\mathbb Q}} \renewcommand{\P}{{\mathbb P}}
\newcommand {\PP}{{\mathbb P}} \newcommand {\CC}{{\mathbb C}}

\newcommand {\kk}{{\mathbf k}}
\newcommand {\mm}{{\boldsymbol \mu}}

\newcommand{\oooo}{\multiput(0,0)(10,0){4}{\circle*{2}}}

\newcommand{\Eeee}{\put(1,0){\line(1,0){8}}}
\newcommand{\eEee}{\put(11,0){\line(1,0){8}}}
\newcommand{\eeEe}{\put(21,0){\line(1,0){8}}}

\DeclareMathOperator {\ev}{ev}

\DeclareMathOperator {\val}{val}

\DeclareMathOperator {\trop}{trop}

\DeclareMathOperator {\mult}{mult}
\DeclareMathOperator {\dive}{div}
\DeclareMathOperator {\floor}{floor}
\DeclareMathOperator {\rel}{rel}
\DeclareMathOperator {\pt}{pt}
\DeclareMathOperator {\codim}{codim}

\renewenvironment {enumerate}%
  {\rule{1mm}{0mm}\begin {oldenumerate}%
    \parskip0.5ex plus0.2ex \itemsep 0mm \parindent 0mm}%
  {\end {oldenumerate}}

\renewenvironment {itemize}%
  {\rule{1mm}{0mm}\begin {olditemize}%
    \parskip0.5ex plus0.2ex \itemsep 0mm \parindent 0mm}%
  {\end {olditemize}}

%% file: introduction.tex
\section {Introduction}

On the moduli spaces $\overline{M}_{g,r}$ and $ \overline{M}_{g,r} (\P^s,d) $
of $r$-marked genus-$g$ stable curves (resp.\ stable maps of degree $d$ to
projective space $ \P^s $), the \emph{Psi-class} $\psi_i$ for $ i=1,\dots,r $
is the first Chern class of the line bundle whose fiber over a point
$(C,x_1,\ldots,x_r)$ (resp.\ $(C,x_1,\ldots,x_r,f)$) is the cotangent space of
$C$ at $x_i$. These Psi-classes are useful to count curves with tangency
conditions, for example. To count curves that satisfy incidence conditions
(e.g.\ pass through given points), one defines evaluation maps on the space of
stable maps, $\ev_i: \overline M_{g,r}(\P^s,d)\rightarrow \P^s$, which send a
stable map $(C,x_1,\ldots,x_r,f)$ to the image $f(x_i)$ of the marked point
$x_i$. Then we can pull back the incidence conditions via the evaluation maps.
Finally, we can intersect pullbacks along the evaluation maps and Psi-classes
on $ \overline M_{g,r}(\P^s,d)$. The degrees of such zero-dimensional
intersection products are called \emph {descendant Gromov-Witten invariants}.
They have been studied in detail in Gromov-Witten theory.

Tropical geometry has been applied to enumerative problems very successfully.
Grigory Mikhalkin has pioneered the field in \cite{Mi03}, proving the
Correspondence Theorem for the numbers $N(d,g)$ of degree-$d$ genus-$g$ nodal
plane curves through $3d+g-1$ points in general position: counting such curves
in algebraic geometry and in tropical geometry will give the same results. He
also developed the tropical \emph {lattice path algorithm} to determine these
numbers. Tropical analogues of moduli spaces of stable curves and maps have
been introduced in \cite{Mi07,GKM07}, and tropical intersection theory was used
to define tropical enumerative numbers for rational curves analogously to the
classical world. Andreas Gathmann and Hannah Markwig showed that the famous
recursion formulas for the count of plane curves known as \emph {Kontsevich's
formula} \cite{GM053} resp.\ \emph {the Caporaso-Harris algorithm} \cite{GM052}
also hold in the tropical world and can be proven using purely tropical
methods. Tropical analogues of Psi-classes on the space of abstract tropical
curves $\mathcal{M}_{0,r}$ have been introduced by Grigory Mikhalkin
\cite{Mi07}, and tropical descendant Gromov-Witten invariants on
$\mathcal{M}_{0,r}(\RR^2,d)$ by Hannah Markwig and Johannes Rau \cite{MR08}.
Markwig and Rau show that these tropical descendant Gromov-Witten invariants
for which every Psi-condition $\psi_i$ comes together with a point condition $
\ev_i^* \pt $ satisfy the so-called \emph {WDVV equations} which can be thought
of as generalizations of Kontsevich's formula. It follows that those numbers
are equal to their classical counterparts, i.e.\ a correspondence theorem holds
here as well (proved indirectly). Tropical curves contributing to the count of
such descendant Gromov-Witten invariants have higher-valent vertices at the
marked points satisfying the Psi-conditions. Markwig and Rau also generalized
the lattice path algorithm to count tropical descendant Gromov-Witten
invariants.

Tropical descendant Gromov-Witten invariants play an important role in a recent
work of Mark Gross \cite{Gro09}: he proves a correspondence theorem for certain
tropical descendant invariants and period integrals on the mirror of $\P^2$.
The philosophy of his paper is that mirror symmetry (or, more precisely, the
correspondence of Gromov-Witten invariants and period integrals on the mirror)
should follow easily from tropical geometry by proving correspondence theorems
for tropical Gromov-Witten invariants and classical Gromov-Witten invariants on
the one hand, and correspondence theorems for tropical Gromov-Witten invariants
and period integrals on the other hand.

The Caporaso-Harris algorithm counts plane curves satisfying point conditions
and multiplicity conditions to a fixed line, resulting in the so-called \emph
{relative Gromov-Witten invariants} of the plane. The rough idea of the
algorithm is to move one of the points from its general position to the line.
After this, the points are no longer in general position, and the curves
satisfying the conditions might split into several components. One then
collects the contributions from all the components and thus produces recursive
relations. An analogue of the Caporaso-Harris algorithm for rational descendant
Gromov-Witten invariants has been developed by Andreas Gathmann \cite {Gat02}.

To apply the same strategy in tropical geometry, we choose the infinitely far
left vertical line. Tropical curves with higher multiplicities to this line are
then just curves with \emph{thick ends}, i.e.\ ends of higher weight in
direction $(-1,0)$. Instead of moving a point to the line, we just move the
point far away to the left. The new point configuration is still in tropically
general position, and the curves satisfying the conditions do not split into
several components. However, if the far left point is not on an end, the
tropical curve contains a part on the far left called a \emph{floor} with one
end in direction $(0,-1)$ and one end in direction $(1,1)$, and this part is
connected to the rest of the curve by horizontal edges only. The
Caporaso-Harris recursion also holds for lattice paths \cite{GM052}.

By applying the Caporaso-Harris algorithm in tropical geometry several times
(i.e.\ spreading the points $ p_1,\dots,p_r $ such that $p_{i+1}$ is far left
of $p_i$ for all $i$) we can decompose the tropical curve into floors. The data
of a tropical curve satisfying these point conditions can then be compressed
into a \emph {floor diagram} as introduced by Erwan Brugall\'e and Grigory
Mikhalkin \cite{BM1,BM2} and studied further by Sergey Fomin and Grigory
Mikhalkin \cite{FM}. They also introduced floor diagrams to count relative
plane Gromov-Witten invariants. Using floor diagrams, Fomin and Mikhalkin have
been able to prove new results about node polynomials, and Florian Block
computed node polynomials for curves with up to $14$ nodes \cite{FB}. Floor
diagrams have also been applied to deduce recursive formulas of Caporaso-Harris
type for Welschinger invariants \cite{ABM08}.

The aim of this paper is to introduce floor diagrams for plane descendant
Gromov-Witten invariants (such that every Psi-condition $\psi_i$ comes together
with a point condition $ \ev_i^* \pt $) which we call \emph {Psi-floor
diagrams}. The count of these diagrams gives exactly the tropical descendant
Gromov-Witten invariants. Because of the Correspondence Theorem it then follows
that they also give the classical descendant Gromov-Witten invariants. We
generalize our definition to relative Psi-floor diagrams and prove that their
count computes tropical relative Gromov-Witten invariants. Afterwards, we show
that the numbers of relative Psi-floor diagrams satisfy a Caporaso-Harris
formula and we show that our formula coincides with the classical formula by
Gathmann mentioned above. It follows that relative Psi-floor diagrams (and thus
also tropical relative descendant Gromov-Witten invariants) count relative
descendant Gromov-Witten invariants.

The difficulty in generalizing the definition of floor diagrams to tropical
curves satisfying Psi-conditions is that, because of the higher-valent
vertices, we cannot necessarily split the curve into single floors. So we have
to introduce multiple floors which are harder to deal with combinatorially. As
a consequence, there is no longer a bijection between labeled floor diagrams
and tropical curves. Rather, there are several tropical curves encoded in one
Psi-floor diagram since there are many ways how a multiple floor can look in a
tropical curve. Thus, we have to introduce new multiplicities for Psi-floor
diagrams that encode how many tropical curves correspond to one diagram.

One can think of tropical geometry as a degeneration of classical geometry, and
it is remarkable that enumerative numbers survive this degeneration. By passing
from tropical curves to floor diagrams, we degenerate even further keeping only
the combinatorial essence of the tropical curve count. Still, this data is
enough to recover the Caporaso-Harris formula. We hope that Psi-floor diagrams
will be useful in the future to prove new results about plane descendant
Gromov-Witten invariants. 

This paper is organized as follows: in Section \ref {sec-classnumbers} we
recall the algebro-geometric definition of absolute and relative descendant
Gromov-Witten invariants and the Caporaso-Harris formula for relative
descendant Gromov-Witten invariants. Correspondingly, we then recall the
definition of tropical descendant Gromov-Witten invariants and their equality
to the corresponding classical numbers in Section \ref {sec-tropnumbers}. We
also generalize this definition to tropical relative descendant Gromov-Witten
invariants. In Section \ref {sec-floor} we introduce Psi-floor diagrams and
their relative analogues and prove that they count the corresponding tropical
curves. We prove that Psi-floor diagrams satisfy the same Caporaso-Harris
formula as the corresponding relative descendant Gromov-Witten invariants. It
follows that relative Psi-floor diagrams (and thus, tropical relative
descendant Gromov-Witten invariants) count relative descendant Gromov-Witten
invariants.

Part of this work was accomplished at the Mathematical Sciences Research
Institute (MSRI) in Berkeley, CA, USA, during the one-semester program on
tropical geometry. The authors would like to thank the MSRI for hospitality. In
addition, Florian Block was supported by the NSF grant DMS-055588, Andreas
Gathmann by the Simons Professorship of the MSRI, and Hannah Markwig by the
MSRI and the German Research Foundation (Deutsche Forschungsgemeinschaft)
through the Institutional Strategy of the University of G\"ottingen.

%% file: classnumbers.tex
\section {Descendant Gromov-Witten invariants} \label {sec-classnumbers}

Let us start by recalling the algebro-geometric construction and computation of
the absolute and relative descendant Gromov-Witten invariants whose
corresponding tropical version we will study later in this paper. For details
in this section we refer mainly to \cite {FP97,KM98} in the absolute and \cite
{Gat02} in the relative case. Throughout this section we will work with the
ground field $ \CC $ of the complex numbers and denote by $ A_*(X) $ and $
A^*(X) $ the Chow homology and cohomology groups of a scheme (or stack) $X$. A
class $ \gamma \in A^i (X) $ will be said to have codimension $ \codim \gamma =
i $, and the class of a hyperplane in a projective space $ \PP^s $ will be
denoted $ h \in A^1(\PP^s) $.

\subsection {Absolute descendant Gromov-Witten invariants}
  \label {subsec-classnumbers-abs}

For $ s>0 $ and $ r,d \ge 0 $ we denote by $ \bar M_{0,r} (\PP^s,d) $ the
moduli space of $r$-marked rational stable maps of degree $d$ to the projective
space $ \PP^s $ (see \cite {FP97} Section 4). Its points correspond to tuples
$ (C,x_1,\dots,x_r,f) $ (modulo automorphisms) such that
\begin {itemize}
\item $C$ is a connected, complete rational curve with at most nodes as
  singularities;
\item $ x_1,\dots,x_r $ are distinct smooth points on $C$;
\item $ f: C \to \PP^s $ is a morphism of degree $d$, i.e.\ such that $ f_* [C]
  $ is the class of $d$ times a line; and
\item the tuple $ (C,x_1,\dots,x_r,f) $ has only finitely many automorphisms.
\end {itemize}
Intuitively, $ \bar M_{0,r} (\PP^s,d) $ can be thought of as a compactification
of the space of all rational degree-$d$ curves in $ \PP^s $ with $r$ marked
points. It is a smooth, complete, and separated stack of dimension $ (s+1)d
+s-3+r $.

For $ i=1,\dots,r $ there are so-called evaluation maps $ \ev_i: \bar M_{0,r}
(\PP^s,d) \to \PP^s $ that send a tuple $ (C,x_1,\dots,x_r,f) $ to the image $
f(x_i) $ of the $i$-th marked point. Moreover, we denote by $ \psi_i \in A^1
(\bar M_{0,r}(\PP^s,d)) $ the $i$-th cotangent line class (also called the
$i$-th Psi-class), i.e.\ the first Chern class of the line bundle whose fiber
over a point $ (C,x_1,\dots,x_r,f) $ is the cotangent space of $C$ at the
(smooth) point $ x_i $.

In general, descendant Gromov-Witten invariants are now defined by taking
degrees of zero-dimensional intersection products of Psi-classes and pull-backs
by the evaluation maps on the above moduli spaces. More precisely, pick $ a_1,
\dots,a_r \ge 0 $ and $ \gamma_1,\dots,\gamma_r \in A^*(\PP^s) $ such that
the dimension condition
  \[ \sum_{i=1}^r (a_i + \codim \gamma_i) = \dim \bar M_{0,r} (\PP^s,d) \]
holds. Then we define the corresponding Gromov-Witten invariant
  \[ \langle \tau^{a_1}(\gamma_1) \cdots \tau^{a_r}(\gamma_r) \rangle_d^{\PP^s}
       := \deg \big( \ev_1^* \gamma_1 \cdot \psi_1^{a_1}
       \cdot \; \cdots \; \cdot
       \ev_r^* \gamma_r \cdot \psi_r^{a_r} \cdot
       [\bar M_{0,r} (\PP^s,d)] \big) \in \QQ. \]
For $ a_1=\cdots=a_r=0 $ we can simply think of this invariant as the number of
rational degree-$d$ curves in $ \PP^s $ passing through $r$ given generic
subvarieties of classes $ \gamma_1,\dots,\gamma_r $. For other choices of
$ a_1,\dots,a_r $ these numbers do not have an immediate geometric
interpretation, but they occur e.g.\ in the computation of numbers of curves
satisfying tangency conditions in addition to incidence conditions. 

The Gromov-Witten invariants above are all well-known; they can be computed
e.g.\ using the WDVV and topological recursion relations (see \cite {KM94}
Section 3, \cite {KM98} Corollary 1.3). In what follows we will need in
particular the following invariants of $ \PP^1 $.

\begin {lemma} \label {lem-inv-p1}
  For all $ a,b,c,d \ge 0 $ with $ a = 2d-2+b $ we have
    \[ \langle \underbrace {1 \, \cdots \, 1}_b \,
               \underbrace {h \, \cdots \, h}_c \, \tau^a(h)
       \rangle_d^{\PP^1} = \frac {d^c}{d!^2}. \]
\end {lemma}

\begin {proof}
  The equation $ a = 2d-2+b $ is simply the dimension condition. Let us first
  assume that $ d>0 $. By the fundamental class and divisor axioms of
  Gromov-Witten invariants (see e.g.\ \cite {Get98} Proposition 12) we then
  know that
    \[ \langle \underbrace {1 \, \cdots \, 1}_b \,
               \underbrace {h \, \cdots \, h}_c \, \tau^a(h)
       \rangle_d^{\PP^1} = d^c \cdot \langle \tau^{a-b}(h) \rangle_d^{\PP^1}.
    \]
  As the one-point invariant $ \langle \tau^a(h) \rangle_d^{\PP^1} $ is equal
  to $ \frac 1{d!^2} $ by \cite {Pan98} Section 1.4, the result follows.

  In the special case $ d=0 $ we see first of all that we must have $ b \ge 2 $
  by the dimension condition. Thus we can again use the fundamental class and
  divisor axioms to reduce the invariant to
    \[ \langle \underbrace {1 \, \cdots \, 1}_b \,
               \underbrace {h \, \cdots \, h}_c \, \tau^a(h)
       \rangle_d^{\PP^1} =
       d^c \cdot \langle 1 \; 1 \; \tau^{0}(h) \rangle_d^{\PP^1} = d^c
    \]
  as stated in the lemma (i.e.\ to $1$ for $ c=0 $ and to $0$ otherwise).
\end {proof}

Our main concern in this paper, however, will be the Gromov-Witten invariants
of the projective plane $ \PP^2 $ where each of the classes $
\gamma_1,\dots,\gamma_r $ above is the class $ \pt = h^2 $ of a point. By the
dimension condition we then need non-negative integers $ a_1,\dots,a_r $ such
that
  \[ 2r+a_1+\dots+a_r = \dim \bar M_{0,r} (\PP^2,d), \qquad
     \mbox {i.e.} \;\; a_1+\dots+a_r = 3d-1-r \]
to get a well-defined number $ \langle \tau^{a_1} (\pt) \, \cdots \, \tau^{a_r}
(\pt) \rangle_d^{\PP^2} $. Note that by the symmetry of the points this number
depends only on how often each Psi-power occurs among the numbers $
a_1,\dots,a_r $. Let us therefore introduce a simplified notation that reflects
this symmetry and that will be particularly useful when considering floor
diagrams later:

\begin {notation}[Sequences] \label {not-sequences}
  Let $ \kk = (\kk_0, \kk_1, \kk_2, \dots)$ be a sequence of non-negative
  integers with only finitely many non-zero entries. We set
  \begin {align*}
    |\kk| &:= \kk_0 + \kk_1 + \kk_2 + \cdots, \\
    I\kk  &:= 0 \, \kk_0 + 1 \, \kk_1 + 2 \, \kk_2 + \cdots, \\
    I^\kk &:= 0^{\kk_0} \cdot 1^{\kk_1} \cdot 2^{\kk_2} \cdot \; \cdots, \\
    \kk!  &:= \kk_0! \cdot \kk_1! \cdot \kk_2! \cdot \; \cdots.
  \end {align*}
  Moreover, if $ \kk,\kk' $ are two such sequences we define the sequence
  $ \kk+\kk' $ by componentwise addition and write $ \kk \le \kk' $ if $ \kk_i
  \le \kk'_i $ for all $ i \ge 0 $. To simplify notation, we will usually write
  such sequences as \emph {finite} sequences $ (\kk_0,\dots,\kk_n) $ for some
  $n$ with the convention that the remaining entries $ \kk_{n+1}, \kk_{n+2},
  \dots $ are then equal to zero.
\end {notation}

\begin {definition}[$ \tilde N_{d,\kk} $ and $ N_{d,\kk} $] \label {def-ndk}
  Let $ d \ge 0 $, and let $ \kk=(\kk_0,\kk_1,\kk_2,\dots) $ be a sequence of
  non-negative integers such that $ I\kk = 3d-1-|\kk| $. For $ r=|\kk| $ let $
  a_1,\dots,a_r $ be an $r$-tuple of non-negative integers that contains each
  number $ i \in \NN $ exactly $ \kk_i $ times (in any order). We then define
    \[ \tilde N_{d,\kk} := \langle \tau^{a_1} (\pt) \, \cdots \, \tau^{a_r}
         (\pt) \rangle_d^{\PP^2}
       \qquad \text {and} \qquad
       N_{d,\kk} := \frac {|\kk| \, !}{\kk!} \, \tilde N_{d,\kk}. \]
\end {definition}

\begin {remark} \label {rem-ndk}
  The two sets of numbers $ \tilde N_{d,\kk} $ and $ N_{d,\kk} $ count the
  following: the $\tilde N_{d,\kk}$ are the numbers of rational plane
  degree-$d$ curves passing through $r$ points and satisfying in addition a $
  \psi^i $ condition at $ \kk_i $ chosen marked points for all $i$. If we do
  not choose the points for the $ \psi^i $ conditions but rather only require
  that among the $r$ marked points there are $ \kk_i $ of them at which a $
  \psi^i $ condition is satisfied (i.e.\ sum over all tuples $ a_1,\dots,a_r $
  above containing each $i$ exactly $ \kk_i $ times) then we get instead the
  numbers $ N_{d,\kk} $, which will turn out to be more natural when
  considering floor diagrams later.
\end {remark}

\subsection {Relative descendant Gromov-Witten invariants}
  \label {subsec-classnumbers-rel}

Relative invariants are very similar to the absolute invariants of Section
\ref {subsec-classnumbers-abs}, except that we now fix once and for all a line
$ H \subset \PP^2 $ and count curves in $ \PP^2 $ that have prescribed local
intersection multiplicities with $H$ in addition to satisfying the evaluation
and Psi-conditions above.

More precisely, choose $ d>0 $ and let $ \mu_1,\dots,\mu_r \in \NN $ for some $
r>0 $ such that $ \mu_1+\cdots+\mu_r = d $. Setting $ \mm = (\mu_1,\dots,
\mu_r) $, we denote by $ \bar M_{0,\mm} (\PP^2,d) \subset \bar M_{0,r}
(\PP^2,d) $ the closure of the subset of all $ (C,x_1,\dots,x_r,f) $ such that
$C$ is smooth and $ f^* H = \mu_1 x_1 + \cdots + \mu_r x_r $ as divisors on $C$
(see \cite {Gat02} Section 1). These spaces are called the moduli spaces of
stable maps relative to $H$; they have dimension $ 2d-1+r $.

As in the absolute case, degrees of zero-dimensional intersection products of
Psi-classes and pull-backs by the evaluation maps on the moduli spaces of
relative stable maps are called relative descendant Gromov-Witten invariants.
So if we now fix $ a_1,\dots,a_r \in \NN $ and $ \gamma_1,\dots,\gamma_r \in
A^* (\PP^2) $ such that
  \[ \sum_{i=1}^r (a_i + \codim \gamma_i) = \dim \bar M_{0,\mm} (\PP^2,d), \]
we can define in a similar way as above an associated relative Gromov-Witten
invariant
  \[ \langle \tau^{a_1}(\gamma_1) \cdots \tau^{a_r}(\gamma_r)
         \rangle_\mm^{\PP^2}
       := \deg (\ev_1^* \gamma_1 \cdot \psi_1^{a_1}
       \cdot \; \cdots \; \cdot
       \ev_r^* \gamma_r \cdot \psi_r^{a_r} \cdot
       [\bar M_{0,\mm} (\PP^2,d)]) \in \QQ. \]
If $ a_1=\cdots=a_r=0 $ this invariant can be interpreted by construction as
the number of plane rational degree-$d$ curves with $r$ marked points that have
local intersection multiplicity $ \mu_i $ and, in addition, pass through a
generic subvariety of $ \PP^2 $ of class $ \gamma_i $ at the $i$-th marked
point, for all $ i=1,\dots,r $. In particular, the marked points $ x_i $ with $
\mu_i>0 $ will lie on $H$, whereas the ones with $ \mu_i=0 $ in general do not.

As before, we will restrict our attention in this paper to a certain subset of
these invariants. Namely, we will only consider choices of $ \mu_1,\dots,\mu_r,
a_1,\dots,a_r,\gamma_1,\dots,\gamma_r $ corresponding to Psi-conditions only at
points away from $H$, i.e.\ such that for all $ i=1,\dots,r $ we have one of
the following cases:
\begin {itemize}
\item $ \mu_i>0 $, $ a_i=0 $, and $ \gamma_i=h^1 $ the class of a line (i.e.\ a
  marked point lying on a fixed point of $H$ with a given local intersection
  multiplicity of the curve to $H$). For $ j \ge 1 $ we will denote the number
  of such $i$ with $ \mu_i = j $ by $ \alpha_j $.
\item $ \mu_i>0 $, $ a_i=0 $, and $ \gamma_i=h^0 $ (i.e.\ a marked point lying
  on a non-fixed point of $H$ with a given local intersection multiplicity of
  the curve to $H$). For $ j \ge 1 $ we will denote the number of such $i$ with
  $ \mu_i = j $ by $ \beta_j $.
\item $ \mu_i=0 $ and $ \gamma_i=h^2 $ (i.e.\ a marked point lying on a fixed
  generic point of $ \PP^2 $ and possibly satisfying some Psi-conditions).
  For $ j \ge 0 $ we will denote the number of such $i$ with $ a_i = j $ by
  $ \kk_j $.
\end {itemize}
By symmetry of the marked points, the three sequences $ \alpha = (\alpha_1,
\alpha_2,\dots) $, $ \beta = (\beta_1,\beta_2,\dots) $, and $ \kk = (\kk_0,
\kk_1,\kk_2,\dots) $ determine the invariant under consideration uniquely.
So we can make the following definition:

\begin {definition}[$ \tilde N_{d,\kk}(\alpha,\beta) $ and $ N_{d,\kk}
    (\alpha,\beta) $] \label {def-ndk-ab}
  With notations as above, we set
    \[ \tilde N_{d,\kk} (\alpha,\beta) :=
       \langle \tau^{a_1}(\gamma_1) \cdots \tau^{a_r}(\gamma_r)
       \rangle_\mm^{\PP^2}. \]
  So $ \tilde N_{d,\kk} (\alpha,\beta) $ is the number of plane rational
  marked degree-$d$ curves $ (C,x_1,\dots,x_r,f) $ satisfying the following
  conditions:
  \begin {itemize}
  \item For each $ i \in \NN $ fix $ \alpha_i $ of the marked points on $C$ and
    a general point on $H$ for each of them; each of these marked points then
    has to be mapped by $f$ to the corresponding given point on $H$, and $C$
    must have local intersection multiplicity $i$ to $H$ there.
  \item For each $ i \in \NN $ fix $ \beta_i $ of the marked points on $C$;
    each of these marked points then has to be mapped by $f$ to $H$, and $C$
    must have local intersection multiplicity $i$ to $H$ there.
  \item For each $ i \in \NN $ fix $ \kk_i $ of the marked points on $C$ and
    a general point in $ \PP^2 $ for each of them; each of these marked points
    then has to be mapped by $f$ to the corresponding given point in $ \PP^2 $,
    and $C$ must satisfy in addition a $ \psi^i $ condition there.
  \end {itemize}
  Note that the dimension condition translates to
    \[ I (\alpha + \beta + \kk) = 3d-1+|\beta|-|\kk| \]
  in these variables, where we use notation \ref {not-sequences} also for
  the sequences $ \alpha $ and $ \beta $ (although they start at index $1$
  rather than $0$). In the same way, the condition $ \mu_1+\cdots+\mu_r = d $
  translates to
    \[ I(\alpha + \beta) = d. \]

  As in Definition \ref {def-ndk} let us also introduce a slight variant of
  these invariants where we do not specify which Psi-power condition has to be
  satisfied at which point $ x_i $ with $ \mu_i = 0 $, and where we do not
  mark the non-fixed points on $H$ of the curves: we set
    \[ N_{d,\kk} (\alpha,\beta)
         := \frac 1{\beta!} \cdot \frac {|\kk|\,!}{\kk!} \cdot
            \tilde N_{d,\kk} (\alpha,\beta). \]
\end {definition}

Just like their absolute counterparts all relative Gromov-Witten invariants
that we have introduced in this section are actually known to be computable
recursively. To do so one uses a generalization of the Caporaso-Harris formula
of \cite {CH98} that we will describe now.

\subsection {The Caporaso-Harris formula for descendant invariants}
  \label {subsec-classnumbers-ch}

In this section we want to use relative Gromov-Witten theory to derive a
recursive formula for the numbers $ \tilde N_{d,\kk} (\alpha,\beta) $ (and thus
also for $ N_{d,\kk} (\alpha,\beta) $) of Definition \ref {def-ndk-ab}.

As in the beginning of Section \ref {subsec-classnumbers-rel} let $ r,d>0 $
and $ \mu_1,\dots,\mu_r \ge 0 $ with $ \mu_1+\cdots+\mu_r=d $. We have then
constructed a moduli space $ \bar M_{0,\mm} (\PP^2,d) \subset \bar M_{0,r}
(\PP^2,d) $ of dimension $ 2d-1+r $ of plane rational degree-$d$ stable maps
relative to a fixed line $ H \subset \PP^2 $, and our invariants $ \tilde
N_{d,\kk} (\alpha,\beta) $ were certain zero-dimensional intersection products
on these spaces.

Since $ \mu_1+\cdots+\mu_r=d $ there can be at most $d$ marked points $ x_i $
with $ \mu_i>0 $. Note that our invariants had no Psi-conditions and at most a
codimension-1 evaluation condition at all these points. So the conditions at
these marked points yield a cycle of codimension at most $d$ --- and as the
dimension of our moduli space is $ 2d-1+r>d $ it follows that there must be at
least one marked point $ x_i $ with $ \mu_i=0 $. By symmetry we may assume
without loss of generality that $ x_1 $ is such a marked point, i.e.\ that $
\mu_1 = 0 $. For our invariant this marked point $ x_1 $ is then required to
map to a given general point in $ \PP^2 $.

The idea of the proof is now to move this generic chosen point to a special
position, namely to a point on $H$. As we have marked all intersection points
of the curves with $H$ already (note that $ \mu_1+\cdots+\mu_r=d $) this
forces the curves to become reducible and split up into several components of
smaller degree, one of which will be mapped completely to $H$. The curves can
then be enumerated recursively over the degree.

To describe this process more formally we follow the notation and results from
Section 2 of \cite {Gat02}. Note, however, that our current situation is a
little simplified compared to \cite {Gat02} since we have assumed here that
$ \mu_1+\cdots+\mu_r = d $.

\begin {construction}[Moduli spaces $ D(A,B) $, see Definition 2.2 of \cite
    {Gat02}] \label {constr-divisor}
  Fix $ r,d>0 $ and a moduli space $ \bar M_{0,\mm}(\PP^2,d) \subset
  \bar M_{0,r} (\PP^2,d) $ with $ \mm = (\mu_1,\dots,\mu_r) $ and $
  \mu_1+\cdots+\mu_r = d $ as above.

  Choose a partition $ A=(A',A^1,\dots,A^t) $ of $ \{1,\dots,r\} $ for some $ t
  \ge 0 $, and let $ \mm^i $ for $ i=1,\dots,t $ be the tuple of all $ \mu_j $
  with $ j \in A^i $ (in any order). Moreover, pick a $ (t+1) $-tuple $
  B=(d',d^1,\dots,d^t) $ of non-negative integers with $ d^i>0 $ for $
  i=1,\dots,t $ and $ d'+d^1+\cdots+d^t = d $. We assume that we have made our
  choices so that
  \begin {equation} \label {eq-divisor-a}
    m^i := d^i - \sum_{j \in A^i} \mu_j > 0
  \end {equation}
  for all $ i=1,\dots,t $, and thus (by adding all these equations up and
  comparing the sum to $ \mu_1+\cdots+\mu_r=d $) so that
  \begin {equation} \label {eq-divisor-b}
    d' + m^1 + \cdots + m^t = \sum_{j \in A'} \mu_j.
  \end {equation}
  In this case we now define the space $ D(A,B) $ to be
    \[ D(A,B) := \bar M_{0,t+\#A'} (H,d') \times_{(\PP^1)^t}
         \prod_{i=1}^t \bar M_{0,(m^i) \cup \mm^i} (\PP^2,d^i), \]
  where $ (m^i) \cup \mm^i $ denotes the $ (\#A^i+1) $-tuple obtained by
  prepending $ m^i $ at the beginning of $ \mm^i $, and the maps to $ (\PP^1)^t
  $ for the fiber product are the evaluation at the first $t$ marked points of
  the first factor and at the first marked point of each of the moduli spaces
  in the second factor. Note that the first factor is a moduli space of
  absolute stable maps to the line $ H \cong \PP^1 $, whereas the second factor
  consists of moduli spaces of stable maps to $ \PP^2 $ relative to $H$.

  By construction, $ D(A,B) $ parameterizes stable maps to $ \PP^2 $ with
  (generically) $ t+1 $ irreducible components: one ``central'' component in
  $H$, and $t$ ``external'' components in $ \PP^2 $ all attached to the central
  one at a point where they have a local intersection multiplicity to $H$
  as given by $ m^1,\dots,m^t $. The $ (t+1) $-tuples $A$ and $B$ simply
  parameterize how the marked points and the degree split up onto the $ t+1 $
  components. In this way $ D(A,B) $ can be considered as a closed subspace
  of $ \bar M_{0,r} (\PP^2,d) $.

  Note that the case $ t=0 $ is allowed (i.e.\ there may be no external
  components at all), as well as $ d'<1 $ and $ d'>1 $ (i.e.\ the central
  component may be a contracted one or a multiple cover of $H$). The following
  picture shows an example of a general element $ (C,x_1,\dots,x_5,f) \in
  D(A,B) $ for $ d=5 $, $ r=3 $, $ \mm=(4,0,1) $, $ A=(\{1\},\emptyset,\{2,3\})
  $, $ B=(1,2,2) $, and thus $ \mm^1 = () $, $ \mm^2 = (0,1) $, $ m^1=2 $, and
  $ m^2=1 $.

  \begin {center} \input {pics/relstable} \end {center}
\end {construction}

The importance of these moduli spaces comes from the fact that they describe
precisely the curves appearing when moving a marked point from a general
position in $ \PP^2 $ to $H$. In fact, all $ D(A,B) $ are divisors in $ \bar
M_{0,\mm} (\PP^2,d) $, and we have the following statement:

\begin {proposition}[Theorem 2.6 of \cite {Gat02}] \label {prop-relstable}
  With notations as above, we have
    \[ \ev_1^* H \cdot \bar M_{0,\mm} (\PP^2,d) =
         \sum_{t,A,B} \frac {m^1 \cdot \; \cdots \; \cdot m^t}{t!} \; D(A,B) \]
  in the Chow group of $ \bar M_{0,\mm} (\PP^2,d) $, where the sum is taken
  over all $ t \ge 0 $, $A$, and $B$ with $ 1 \in A' $ and satisfying
  condition (\ref {eq-divisor-a}) (and thus also (\ref {eq-divisor-b})) as in
  Construction \ref {constr-divisor}.
\end {proposition}

As usual in Gromov-Witten theory it is now convenient to replace the fiber
product in the Construction \ref {constr-divisor} of $ D(A,B) $ by the
``diagonal splitting'' trick: the fiber product $ X \times_{\PP^1} Y $ of two
spaces $X$ and $Y$ with projections $p$ and $q$ to $ \PP^1 $ can be rewritten
as the pull-back of the diagonal of $ \PP^1 \times \PP^1 $ by the map $ p
\times q $, and as this diagonal has class $ h \times 1 + 1 \times h $ it
follows that
  \[ X \times_{\PP^1} Y = (p^* h + q^* h) \cdot (X \times Y). \]
Let us apply this formula in the expression for $ D(A,B) $ from Construction
\ref {constr-divisor} for each of the $t$ factors $ \PP^1 $ over which we take
the fiber product, thus converting $ D(A,B) $ into a sum of $ 2^t $ terms with
no fiber products. By symmetry, we can then always relabel the external $t$
components so that the ones with the $ ev^* h $ term in the $ \bar M_{0,(m^i)
\cup \mm^i} (\PP^2,d^i) $ factor come first --- if there are $ t' \in
\{0,\dots,t\} $ of these components we then have $ \binom t{t'} $ terms in
the diagonal splitting that become the same after this relabeling. Hence we
can rewrite the formula of Proposition \ref {prop-relstable} in the following
form:
\begin {align*}
  \ev_1^* H \cdot \bar M_{0,\mm} (\PP^2,d)
    &=\! \sum_{t,A,B} \sum_{t'=0}^t \;
       \frac {m^1 \cdot \; \cdots \; \cdot m^t}{t'! \, (t-t')!} \,
       \Big(
         \ev_{t'+1}^* h \cdot \; \cdots \; \cdot \ev_t^* h \cdot
         \bar M_{0,t+\#A'} (H,d')
       \Big) \\
    &  \qquad \;\;
       \times \prod_{i=1}^{t'} \Big(
         \ev_1^*h \cdot \bar M_{0,(m^i) \cup \mm^i} (\PP^2,d^i)
       \Big)
       \times \prod_{i=t'+1}^t \bar M_{0,(m^i) \cup \mm^i} (\PP^2,d^i).
\end {align*}
To get a recursive relation for the invariants
  \[ \tilde N_{d,\kk} (\alpha,\beta) =
       \langle \tau^{a_1}(\gamma_1) \cdots \tau^{a_r}(\gamma_r)
       \rangle_\mm^{\PP^2} \]
of Definition \ref {def-ndk-ab} we now intersect this equation of cycles with
the class
  \[ \ev_1^* h \cdot \psi_1^{a_1} \cdot
     \ev_2^* \gamma_2 \cdot \psi_2^{a_2} \cdot \; \cdots \; \cdot
     \ev_r^* \gamma_r \cdot \psi_r^{a_r} \]
(note that $ \gamma_1 = \pt $ by assumption, and thus the two evaluations
$ \ev_1^* H \cdot \ev_1^* h $ together give the desired condition $ \ev_1^*
\gamma_1 $ at the first point). The left hand side of the equation is then
simply $ \tilde N_{d,\kk} (\alpha,\beta) $. Each summand on the right hand side
is a product of one absolute Gromov-Witten invariant of $ \PP^1 $ and $t$
relative Gromov-Witten invariants of $ \PP^2 $. The invariant of $ \PP^1 $
has the condition $ \ev_1^* h \cdot \psi_1^{a_1} $ at the first marked point,
a condition $ \ev_i^* h $ at all gluing points from the last $ t-t' $ external
components and all $ x_i $ with $ i \in A' $ such that $ \gamma_i = h $, and no
condition at all at the other points. On the other hand, the $t$ relative
invariants of $ \PP^2 $ are again of the type of invariants considered in
Definition \ref {def-ndk-ab}: we can write them as $ N_{d^i,\kk^i} (\alpha^i +
e_{m^i},\beta^i) $ for the first $ t' $ and $ N_{d^i,\kk^i} (\alpha^i,\beta^i
+e_{m^i}) $ for the last $ t-t' $ invariants, where $ \alpha^i,\beta^i,\kk^i $
denote the sequences associated to the marked points $ x_j $ with $ j \in A^i $
according to Definition \ref {def-ndk-ab}. Finally, let us then rewrite the sum
over $A$ as a sum over the corresponding sequences $ \alpha^i $, $ \beta^i $,
$ \kk^i $. If we set
\begin {equation} \label {eq-binoms}
  \alpha' := \alpha - \alpha^1 - \cdots - \alpha^t
    \quad \mbox {and} \quad
  \binom {\alpha}{\alpha^1,\dots,\alpha^t} := \prod_{i \ge 1} \,
    \frac {\alpha_i!}{\alpha^1_i! \cdot \; \cdots \; \cdot \alpha^t_i!
    \cdot \alpha'_i!}
\end {equation}
(and similarly for $ \beta $ and $ \kk $, except that the index of the
sequences starts at $0$ for $ \kk $), then exactly
  \[ \binom {\alpha!}{\alpha^1,\dots,\alpha^t} \cdot
     \binom {\beta!}{\beta^1,\dots,\beta^t} \cdot
     \binom {\kk-e_a}{\kk^1,\dots,\kk^t} \]
choices of partitions of $A$ into $t$ subsets will give rise to the same
invariants. Here, $ e_a $ denotes the sequence with only non-zero entry $1$
in the $a$-th component --- we have to write $ \kk-e_a $ instead of $ \kk $
since the first marked point is fixed to lie on the central component, so there
is no choice here where to put this point. Hence our equation becomes
\begin {align*}
  \tilde N_{d,\kk}(\alpha,\beta)
    &= \sum_{t,t'} \sum_{\alpha,\beta,\kk} \sum_B
       \frac {m^1 \cdot \; \cdots \; \cdot m^t}{t'! \, (t-t')!} \,
       \binom {\alpha}{\alpha^1,\dots,\alpha^t} \,
       \binom {\beta}{\beta^1,\dots,\beta^t} \,
       \binom {\kk-e_a}{\kk^1,\dots,\kk^t} \\[2mm]
    &\qquad\quad
       \cdot \langle
         \underbrace {\vphantom {\big(} 1 \, \cdots \, 1}_{|\beta'|+t'} \,
         \underbrace {\vphantom {\big(} h \, \cdots \, h}_{|\alpha'|+t-t'} \,
         \tau^a (h)
       \rangle_{d'}^{\PP^1} \\[2mm]
    &\qquad\quad
       \cdot \prod_{i=1}^{t'} N_{d^i,\kk^i} (\alpha^i+e_{m^i},\beta^i)
       \cdot \prod_{i=t'+1}^{t} N_{d^i,\kk^i} (\alpha^i,\beta^i+e_{m^i})
\end {align*}
Note that we must have $ \kk^1+\cdots+\kk^t = \kk-e_a $ in each term since
marked points with generic point conditions in $ \PP^2 $ cannot lie in the
central component within $H$. Moreover, each relative invariant in this
expression must of course satisfy the dimension condition
  \[ \begin {array}{rr@{\;}c@{\;}l@{\qquad}l}
     & I((\alpha^i+e_{m^i})+\beta^i+\kk^i)
     &=& 3d^i-1+|\beta^i|-|\kk^i|
     &\mbox {for $ i \le t' $} \\
     \mbox {resp.}
     & I(\alpha^i+(\beta^i+e_{m^i})+\kk^i)
     &=& 3d^i-1+|\beta^i+e_{m^i}|-|\kk^i|
     &\mbox {for $ i > t' $}
     \end {array} \]
of Definition \ref {def-ndk-ab}, as well as condition (\ref {eq-divisor-a})
  \[ m^i = d^i - I(\alpha^i+\beta^i) \]
of Construction \ref {constr-divisor}. We can think of the first of these
equations as determining $ d^i $, and of the second as determining $ m^i $ from
$ \alpha^i $, $ \beta^i $, and $ \kk^i $. Finally, inserting the expression
of Lemma \ref {lem-inv-p1} for the absolute Gromov-Witten invariant of $ \PP^1
$ we get the following result that allows us to compute all numbers $ \tilde
N_{d,\kk} (\alpha,\beta) $ recursively.

\begin {theorem}[Caporaso-Harris formula for the relative descendant
    Gromov-Witten invariants $ \tilde N_{d,\kk}(\alpha,\beta) $]
    \label {thm-CHclass}
  The relative Gromov-Witten invariants $ \tilde N_{d,\kk} $ of Definition
  \ref {def-ndk-ab} satisfy the relations
  \begin {align*}
    \tilde N_{d,\kk}(\alpha,\beta)
      &= \sum
         \frac {m^1 \cdot \; \cdots \; \cdot m^t}{t'! \, (t-t')!} \cdot
         \frac {{d'}^{|\alpha'|+t-t'}}{d'!^2} \,
         \binom {\alpha}{\alpha^1,\dots,\alpha^t} \,
         \binom {\beta}{\beta^1,\dots,\beta^t} \,
         \binom {\kk-e_a}{\kk^1,\dots,\kk^t} \\
      &\qquad\quad
         \cdot \prod_{i=1}^{t'}
           \tilde N_{d^i,\kk^i} (\alpha^i+e_{m^i},\beta^i)
         \cdot \prod_{i=t'+1}^t
           \tilde N_{d^i,\kk^i} (\alpha^i,\beta^i+e_{m^i})
  \end {align*}
  for each $ a \in \NN $ with $ \kk_a>0 $. Here, the sum is taken over all
  $ 0 \le t' \le t $ and all sequences $ \alpha^1,\dots,\alpha^t $,
  $ \beta^1,\dots,\beta^t $, $ \kk^1,\dots,\kk^t $ such that
  \begin {itemize}
  \item $ \alpha' := \alpha - \alpha^1 - \cdots - \alpha^t \ge 0 $,
    $ \beta' := \beta - \beta^1 - \cdots - \beta^t \ge 0 $, and
    $ \kk^1 + \cdots + \kk^t = \kk-e_a $;
  \item $ d^i := \frac 13 (I(\alpha^i+\beta^i+\kk^i+e_{m^i})
    -|\beta^i|+|\kk^i|+1) \in \NN_{>0} $ for $ i = 1,\dots,t' $, and
    $ d^i := \frac 13 (I(\alpha^i+\beta^i+\kk^i+e_{m^i})
    -|\beta^i|+|\kk^i|) \in \NN_{>0} $ for $ i = t'+1,\dots,t $;
  \item $ d' := d-d_1-\cdots-d_t \ge 0 $;
  \item $ m^i := d^i - I(\alpha^i+\beta^i) > 0 $ for all $ i=1,\dots,t $.
  \end {itemize}
\end {theorem}

It is easy to rewrite this formula so that it computes the invariants $
N_{d,\kk} (\alpha,\beta) $ instead of $ \tilde N_{d,\kk} (\alpha,\beta) $:

\begin {corollary}[Caporaso-Harris formula for the relative descendant
    Gromov-Witten invariants $ N_{d,\kk}(\alpha,\beta) $] \label {cor-CHclass}
  The invariants $ N_{d,\kk} (\alpha,\beta) $ of Definition \ref {def-ndk-ab}
  satisfy the relations
  \begin {align*}
    N_{d,\kk}(\alpha,\beta)
      &= \sum_{a:\,\kk_a>0} \sum
         \frac {m^1 \cdot \; \cdots \; \cdot m^t}{t'! \, (t-t')!} \cdot
         \frac {{d'}^{|\alpha'|+t-t'}}{d'!^2} \,
         \binom {\alpha}{\alpha^1,\dots,\alpha^t} \,
         \frac 1{\beta'!} \,
         \binom {|\kk|-1}{|\kk^1|,\dots,|\kk^t|} \\
      &\qquad\quad
         \cdot \prod_{i=1}^{t'}
           N_{d^i,\kk^i} (\alpha^i+e_{m^i},\beta^i)
         \cdot \prod_{i=t'+1}^t
           (\beta^i_{m^i}+1) \,
           N_{d^i,\kk^i} (\alpha^i,\beta^i+e_{m^i})
  \end {align*}
  where the second sum is taken over the same partitions and with the same
  conditions as in Theorem \ref {thm-CHclass}.
\end {corollary}

\begin {proof}
  Inserting the expression of Definition \ref {def-ndk-ab} for the numbers $
   N_{d,\kk} (\alpha,\beta) $ in terms of $ \tilde N_{d,\kk} (\alpha,\beta) $
  into the formula of Theorem \ref {thm-CHclass} gives
  \begin {align*}
    N_{d,\kk}(\alpha,\beta)
      &= \sum
         \frac {m^1 \cdot \; \cdots \; \cdot m^t}{t'! \, (t-t')!} \cdot
         \frac {{d'}^{|\alpha'|+t-t'}}{d'!^2} \,
         \binom {\alpha}{\alpha^1,\dots,\alpha^t} \,
         \frac 1{\beta'!} \,
         \binom {|\kk|-1}{|\kk^1|,\dots,|\kk^t|} \,
         \frac {|\kk|}{\kk_a} \\
      &\qquad\quad
         \cdot \prod_{i=1}^{t'}
           N_{d^i,\kk^i} (\alpha^i+e_{m^i},\beta^i)
         \cdot \prod_{i=t'+1}^t
           (\beta^i_{m^i}+1) \,
           N_{d^i,\kk^i} (\alpha^i,\beta^i+e_{m^i})
  \end {align*}
  for all $a$ with $ \kk_a>0 $. Multiplying these equations with $ \frac
  {\kk_a}{|\kk|} $ and summing them up for all $a$ then gives the desired
  equation since $ \sum_a \frac {\kk_a}{|\kk|} = 1 $.
\end {proof}

%% file: pics/relstable.tex
\begin{picture}(0,0)%
\includegraphics{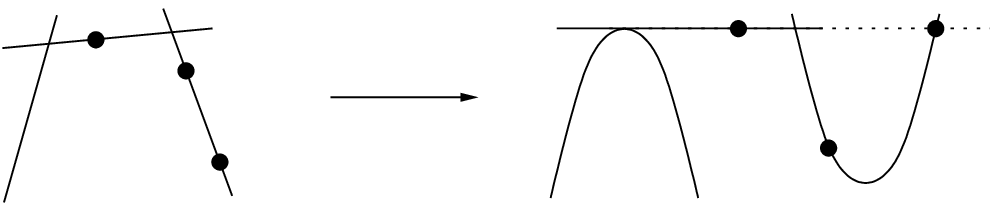}%
\end{picture}%
\setlength{\unitlength}{4144sp}%
\begingroup\makeatletter\ifx\SetFigFont\undefined%
\gdef\SetFigFont#1#2#3#4#5{%
  \reset@font\fontsize{#1}{#2pt}%
  \fontfamily{#3}\fontseries{#4}\fontshape{#5}%
  \selectfont}%
\fi\endgroup%
\begin{picture}(4587,1039)(964,-418)
\put(2791,137){\makebox(0,0)[b]{\smash{{\SetFigFont{10}{12.0}{\rmdefault}{\mddefault}{\updefault}{\color[rgb]{0,0,0}$f$}%
}}}}
\put(1895,142){\makebox(0,0)[lb]{\smash{{\SetFigFont{10}{12.0}{\rmdefault}{\mddefault}{\updefault}{\color[rgb]{0,0,0}$x_2$}%
}}}}
\put(2046,-279){\makebox(0,0)[lb]{\smash{{\SetFigFont{10}{12.0}{\rmdefault}{\mddefault}{\updefault}{\color[rgb]{0,0,0}$x_3$}%
}}}}
\put(1405,423){\makebox(0,0)[b]{\smash{{\SetFigFont{10}{12.0}{\rmdefault}{\mddefault}{\updefault}{\color[rgb]{0,0,0}$x_1$}%
}}}}
\put(4344,474){\makebox(0,0)[b]{\smash{{\SetFigFont{10}{12.0}{\rmdefault}{\mddefault}{\updefault}{\color[rgb]{0,0,0}$f(x_1)$}%
}}}}
\put(4698,-287){\makebox(0,0)[rb]{\smash{{\SetFigFont{10}{12.0}{\rmdefault}{\mddefault}{\updefault}{\color[rgb]{0,0,0}$f(x_2)$}%
}}}}
\put(5217,472){\makebox(0,0)[rb]{\smash{{\SetFigFont{10}{12.0}{\rmdefault}{\mddefault}{\updefault}{\color[rgb]{0,0,0}$f(x_3)$}%
}}}}
\put(5536,344){\makebox(0,0)[lb]{\smash{{\SetFigFont{10}{12.0}{\rmdefault}{\mddefault}{\updefault}{\color[rgb]{0,0,0}$H$}%
}}}}
\put(5536,-331){\makebox(0,0)[rb]{\smash{{\SetFigFont{10}{12.0}{\rmdefault}{\mddefault}{\updefault}{\color[rgb]{0,0,0}$ \PP^2 $}%
}}}}
\put(1486,-331){\makebox(0,0)[b]{\smash{{\SetFigFont{10}{12.0}{\rmdefault}{\mddefault}{\updefault}{\color[rgb]{0,0,0}$C$}%
}}}}
\end{picture}%

%% file: tropnumbers.tex
\section {Tropical descendant Gromov-Witten invariants}
  \label {sec-tropnumbers}

In the last section we have introduced several algebro-geometric descendant
rational Gromov-Witten invariants of the projective plane:
\begin {itemize}
\item the \emph {absolute} invariants $ \tilde N_{d,\kk} $ and $ N_{d,\kk} $
  counting degree-$d$ curves through points and Psi-conditions as specified by
  $ \kk $ (see Definition \ref {def-ndk});
\item the \emph {relative} invariants $ \tilde N_{d,\kk}(\alpha,\beta) $ and
  $ N_{d,\kk}(\alpha,\beta) $ counting degree-$d$ curves through points,
  Psi-conditions as specified by $ \kk $, and multiplicity conditions to a
  fixed line as specified by $ \alpha $ and $ \beta $ (see Definition \ref
  {def-ndk-ab}).
\end {itemize}
The convention here was that the numbers called $ \tilde N $ consider all
points at which some condition has to be satisfied as marked points, whereas
the numbers called $N$ are obtained from these by a simple combinatorial factor
dividing out some symmetries in the conditions.

We will now introduce corresponding numbers with a superscript ``trop'' (e.g.\
$ \tilde N^{\trop}_{d,\kk} $) arising from the count of \emph {tropical
curves}, as well as --- in the following Section \ref {sec-floor} --- numbers
with a superscript ``floor'' (e.g.\ $ \tilde N^{\floor}_{d,\kk}) $ obtained by
counting \emph {floor diagrams}. The convention mentioned above will still hold
for these numbers; we will see however that the $N$ numbers seem to be more
natural from the point of view of floor diagrams, whereas the $ \tilde N $ have
been more natural in the algebro-geometric setting. In the end however, all
corresponding numbers will turn out to be the same, e.g.\ $ \tilde N_{d,\kk} =
\tilde N^{\trop}_{d,\kk} = \tilde N^{\floor}_{d,\kk} $ for all $d$ and $\kk$.
In fact, this is the main result of this paper: that the (rational plane)
absolute and relative descendant Gromov-Witten invariants of algebraic geometry
can also be computed using certain counts of floor diagrams.

\subsection {Absolute tropical descendant Gromov-Witten invariants}

As mentioned in the introduction, tropical descendant Gromov-Witten invariants
can be defined as intersection products on the tropical analogue of the moduli
spaces of stable maps \cite{MR08}. However, in order to avoid introducing too
much notation, we choose to define them here purely in terms of the
combinatorial properties of the tropical curves which we want to count.

A \emph {(rational) abstract tropical curve} is a connected metric graph $
\Gamma $ of genus $0$ (considered as a topological space, with the edges
homeomorphic to closed real intervals), such that unbounded edges (with no
vertex there) are allowed, and such that each vertex has valence at least 3
(see \cite {GKM07} Definition 3.2). The unbounded edges will be called \emph
{ends}, and the length of a bounded edge $e$ will be denoted $ l(e) \in
\RR_{>0} $. We say that such a curve is an \emph {$n$-marked abstract tropical
curve} if $n$ of the ends are marked by $ x_1,\ldots,x_n $. Two (marked)
abstract tropical curves are isomorphic (and will from now on be identified) if
there is an isometry between them (that respects $ x_1,\dots,x_n $ in the
marked case).

We now want to consider maps from marked abstract tropical curves to $ \RR^2 $.
For our later purposes it will be convenient to consider some of the \emph
{left} ends to be marked ends, whereas the other (non-contracted) ends will be
unmarked.

\begin {definition} \label {def-paramcurve}
  Let $ m \ge n \ge 0 $. A \emph {(parameterized plane) $n$-marked tropical
  curve (with $ m-n $ marked left ends)} is a tuple $ (\Gamma,x_1,\ldots,x_m,h)
  $, where $ (\Gamma,x_1,\dots,x_m) $ is an $m$-marked abstract tropical curve
  and $ h: \Gamma \to \RR^2 $ is a continuous map satisfying the following
  conditions.
  \begin {itemize}
  \item On each edge $e$ the map $h$ is integer affine linear, i.e.\ of the
    form $ h(t) = a + t \cdot v $ for $ a \in \RR^2 $ and $ v \in \ZZ^2 $. If
    $ V \in \partial e $ is a vertex of the edge $e$ and we parameterize $e$
    starting at $V$, the vector $v$ in the above equation will be denoted $
    v(V,e) $ and called the \emph {direction vector} of $e$ starting at $V$.
    If $V$ is understood from the context (e.g.\ in case $e$ is an end, having
    only one adjacent vertex) we will also write $ v(e) $ instead of $ v(V,e)
    $. The lattice length of $ v(V,e) $ (i.e.\ the greatest common divisor of
    the entries of $ v(V,e) $) will be called the \emph {weight} $ \omega(e) $
    of $e$.
  \item At each vertex $V$ the balancing condition
      \[ \sum_{e: \, V \in \partial e} v(V,e)=0 \]
    is satisfied.
  \item Each marked end $ x_i $ for $ i=1,\dots,n $ is contracted by $h$ (i.e.\
    $ v(x_i)=0 $).
  \item Each marked end $ x_i $ for $ i=n+1,\ldots,m $ is a \emph {left end}
    (i.e.\ it is of direction $ (-l,0) $ for some $ l\in \NN_{>0} $).
  \end{itemize}
  Two parameterized tropical curves are isomorphic if there is an isomorphism
  of the underlying marked abstract tropical curves commuting with $h$. Note
  that a parametrized $n$-marked curve will in general have more ends than the
  marked ones. The \emph {degree} of such a curve is defined to be the multiset
  consisting of the directions of these non-marked ends, together with the
  directions of the marked left ends $x_{n+1},\dots,x_m$. If the degree
  multiset consists of $d$ copies of each of the vectors $(-1,0)$, $(0,-1)$,
  and $(1,1)$ we say that the curve is \emph {of degree $d$} (see Example \ref
  {ex-nonrelativecurve}).
\end {definition}

\begin {definition}[Multiplicity of a curve] \label {def-mult}
  Let $ C=(\Gamma,x_1,\ldots,x_m,h) $ be a marked tropical curve of degree $
  \Delta = \{v_1,\ldots,v_1,v_2,\ldots,v_2,\ldots,v_r,\ldots,v_r\}$ (with $
  v_1,\dots,v_r $ distinct) such that all vertices that are not adjacent to
  any of the contracted ends $ x_1,\dots,x_n $ are $3$-valent. Let $ V_1,\dots,
  V_t $ be the vertices of $ \Gamma $. For $ i=1,\dots,t $ and $ j=1,\dots,r $
  let $ b_{ij} $ the number of non-marked ends adjacent to $ V_i $ of direction
  $ v_j $.

  As in \cite{Mi03} Definition 2.16 we define the multiplicity of a $3$-valent
  vertex of $C$ to be the absolute value of the determinant of two adjacent
  direction vectors. Setting $ \nu_C := \prod_{i=1}^t \, \prod_{j=1}^r \, \frac
  1{b_{ij}!} $, we then define the \emph {multiplicity} $\mult(C)$ of $C$ to be
  $ \nu_C $ times the product of the multiplicities of all vertices without
  adjacent contracted ends.
\end {definition}

\begin {definition}[$ \tilde N^{\trop}_{d,\kk} $] \label {def-ntroptilde}
  Let $d \ge 1$, and let $\kk$ be a sequence of non-negative integers with $ I
  \kk = 3d -1 - |\kk|$. Furthermore, for $ n = |\kk| $, fix a vector $(a_1,
  \dots,a_n) $ that contains each number $ i \in \NN $ exactly $\kk_i$ times
  (in any order). Let $p_1,\ldots,p_n \in \RR^2$ be points in general position
  (see Definitions 3.2 and 9.7 of \cite{MR08}). We define
    \[ \tilde N^{\trop}_{d,\kk} := \sum_C \mult(C), \]
  where the sum goes over all tropical curves $ C=(\Gamma,x_1,\ldots,x_n,h) $
  (with non-marked left ends, i.e.\ $m-n=0$) of degree $d$ satisfying 
  \begin{itemize}
  \item $h(x_i)=p_i$ for all $i=1,\ldots,n$, and
  \item the end $x_i$ is adjacent to a vertex of valence $a_i+3$ for all $i=1,
    \ldots,n$.
  \end{itemize}
  It follows from the general position of the points that all other vertices of
  $\Gamma$ are then $3$-valent.
\end{definition}

\begin {example} \label {ex-nonrelativecurve}
  The following picture shows a parameterized $9$-marked tropical curve. We
  have drawn the contracted marked ends as dotted lines. We did not specify
  the lengths of the bounded edges in the abstract curve since they are
  determined by the lengths of the images and the (non-zero) direction vectors,
  which in turn are determined by the directions of the ends using the
  balancing condition. The direction vectors are all primitive except for the
  edge with weight $2$ in the image.

  \begin {center} \input {pics/expsicurve} \end {center}

  This curve contributes to $\tilde{N}^{\trop}_{5, {\bf k}}$, where ${\bf k}=
  (7,0,1,1)$, and where we chose $a = (0,0,0,0,0,2,0,3,0)$. Its multiplicity is
  $ \frac{1}{2}\cdot \frac{1}{2}\cdot 2\cdot 2=2$. The two factors of
  $\frac{1}{2}$ arise because two non-marked ends of the same direction are
  adjacent to the end vertex of $x_6$ and of $x_8$. The two factors of $2$ are
  the vertex multiplicities of the vertices of the edges of weight $2$ (not
  adjacent to a contracted end). In the future, we want to avoid drawing the
  abstract curve together with its image. Therefore, we introduce the following
  shortcut for the picture above. When two edges of the abstract curve are
  mapped on top of each other in the image, we choose to draw them separately,
  but close to each other. In this way we can recover the parameterizing
  abstract curve uniquely (see \cite{MR08} Lemma 9.9).

  \begin {center} \input {pics/expsicurve1} \end {center}
\end{example}

For every vector $(a_1, \dots, a_n)$ containing $i$ exactly $\kk_i$ times for
all $i \ge 0$, the number $ \tilde N^{\trop}_{d,\kk}$ equals the tropical
intersection product $ \prod_{i=1}^n \ev_i^*(p_i) \, \psi_i^{k_i} $ on the
moduli space $ \mathcal{M}_{0,n}(\RR^2,d) $ of rational tropical $n$-marked
curves in $ \RR^2 $ of degree $d$ by Remark 3.3 of \cite{MR08}, and is thus a
\emph {tropical descendant Gromov-Witten invariant}.

Later on, it will be convenient to allow arbitrary orderings of the Psi-powers.
This leads to the following invariants.

\begin {definition}[$ N^{\trop}_{d,\kk} $] \label {def-ntrop}
  For $ d\geq 1 $ and $ \kk $ a sequence of non-negative integers with $ I\kk =
  3d-1-|\kk| $ we define the number $ N^{\trop}_{d,\kk} := \sum_C \mult(C) $
  analogously to Definition \ref {def-ntroptilde}, where now the sum is over
  all tropical curves $C$ of degree $d$ with non-marked left ends, such that
  for all $i$ there are $ \kk_i $ contracted ends whose adjacent vertex has
  valence $ i+3 $.

  Obviously, these numbers $ N^{\trop}_{d,\kk}$ are related to the numbers
  $ \tilde N^{\trop}_{d,\kk} $ of Definition \ref {def-ntroptilde} by $
  N^{\trop}_{d,\kk} =  \frac{|\kk|\,!}{\kk!}\tilde{N}^{\trop}_{d,
  \kk} $.
\end {definition}

\begin {remark}[The equality $\tilde N^{\trop}_{d,\kk}=\tilde{N}_{d,\kk}$]
    \label {rem-equiv-trop-alg}
  In \cite{MR08} it was shown that tropical descendant Gromov-Witten invariants
  $\tilde{N}^{\trop}_{d, {\bf k}}$ satisfy the WDVV relations, just as their
  classical counterparts $\tilde{N}_{d, {\bf k}}$ do. As the initial values
  coincide, we can conclude that $\tilde{N}^{\trop}_{d, {\bf k}}=\tilde{N}_{d,
  {\bf k}}$ for all $d$ and $ \kk $. There is no direct bijection of the
  corresponding curves known at this point. Since both pairs of numbers $
  \tilde{N}^{\trop}_{d, {\bf k}}$, $N^{\trop}_{d, {\bf k}}$ and $ \tilde{N}_{d,
  {\bf k}}$, $N_{d, {\bf k}}$ differ by the same combinatorial factor, it
  follows of course that also $N^{\trop}_{d, {\bf k}}=N_{d, {\bf k}}$. Both
  equalities also follow as the special case $ \alpha=() $, $ \beta=(d) $ from
  our Caporaso-Harris formulas (see Remark \ref {rem-eqreltropclass}). 
\end {remark}

\subsection {Relative tropical descendant Gromov-Witten invariants}

For two sequences $\alpha$ and $\beta$ with $d = I(\alpha+\beta) $ let
\begin{multline*}
  \Delta(\alpha,\beta)
    = \{ \underbrace{(-1,0),\ldots,(-1,0)}_{\alpha_1+\beta_1},
         \underbrace{(-2,0),\ldots,(-2,0)}_{\alpha_2+\beta_2},\ldots, \\
         \underbrace{(0,-1),\ldots,(0,-1)}_{d},
         \underbrace{(1,1),\ldots,(1,1)}_{d} \}
\end{multline*}
and consider parameterized $n$-marked tropical curves of degree $ \Delta
(\alpha,\beta)$ with $m-n= |\alpha+\beta|$ marked left ends (i.e.\ all the left ends are marked).

\begin {definition}[$ \tilde N^{\trop}_{d,\kk}(\alpha,\beta) $ and
    $ N^{\trop}_{d,\kk}(\alpha,\beta) $]
  Let $d \ge 1$, and let $\kk$ be a sequence with $I(\alpha+\beta+ \kk) = 3d -
  1+|\beta| - |\kk|$. Furthermore, for $n = |\kk|$ fix a vector $(a_1, \dots,
  a_n)$ containing each $ i \ge 0 $ exactly $\kk_i$ times. Let $p_1,\ldots,p_n
  \in \RR^2$ be points and $ y_{n+1}, \dots,y_{n+|\alpha|}$ be $y$-coordinates
  in general position (analogously to Definitions 3.2 and 9.7 of \cite{MR08}).
  For all $i=n+1,\ldots,n+|\alpha|$ choose a weight $\mu_i$ such that in total
  we have chosen each weight $ k \ge 1 $ exactly $ \alpha_k $ times. In the
  same way, choose weights $ \mu_i $ for $ i = n+|\alpha|+1, \dots,
  n+|\alpha+\beta| $ so that in total we have chosen each weight $ k \ge 1 $
  exactly $ \beta_k $ times.

  We then define
    \[ \tilde N^{\trop}_{d,\kk}(\alpha,\beta)
         := \sum_C \frac 1 {I^\alpha} \, \mult(C), \]
  where the sum is taken over all tropical curves $ C=(\Gamma,x_1,\ldots,x_m,h)
  $ with $m-n= |\alpha+\beta|$ marked left ends (i.e.\ all left ends are
  marked) of degree $ \Delta(\alpha,\beta) $ satisfying
  \begin{itemize}
  \item $h(x_i)=p_i$ for all $i=1,\ldots,n$;
  \item the end $x_i$ is adjacent to a vertex of valence $a_i+3$ for all $i=1,
    \ldots,n$;
  \item for $ i=n+1,\dots,n+|\alpha| $, the $y$-coordinate of $h(x_i)$ equals
    $y_i$;
  \item for $ i=n+1,\dots,n+|\alpha+\beta| $, the marked end $x_i$ is of weight
    $ \mu_i $, i.e.\ we have $ v(x_i) = (-\mu_i,0) $.
  \end{itemize}
  Again, it follows from the general position of the points that all other
  vertices of $\Gamma$ are $3$-valent.

  We also define the numbers $N^{\trop}_{d, {\bf k}}(\alpha, \beta)$
  analogously to Definition \ref{def-ntrop} as numbers of tropical curves
  passing through the given points, with $\kk_i$ contracted ends whose adjacent
  vertex has valence $ i+3 $ for all $i$, with non-marked left ends of the
  specified weights, and satisfying that the prescribed set of $y$-coordinates
  for a given weight are the $y$-coordinates of left ends of this weight. The
  curves are counted with multiplicity $ \frac 1{I^\alpha} \, \mult(C) $ as
  above. The numbers $ N^{\trop}_{d, {\bf k}}(\alpha, \beta)$ and $ \tilde
  N^{\trop}_{d, {\bf k}}(\alpha,\beta)$ are related by $ N^{\trop}_{d,\kk}
  (\alpha,\beta)= \frac{1}{\beta!} \frac{|\kk|}{\kk!}\cdot
  \tilde{N}^{\trop}_{d, {\bf k}}( \alpha,\beta) $.
\end{definition}

Even though tropical descendant Gromov-Witten invariants are defined in
\cite{MR08} only in the non-relative case, a completely analogous argument
shows that the numbers $\tilde{N}^{\trop}_{d, {\bf k}}(\alpha, \beta)$ can also
be interpreted as intersection products of evaluation pull-backs and
Psi-classes on a suitable moduli space of tropical curves. Hence we can think
of these numbers as \emph {tropical relative descendant Gromov-Witten
invariants}.

\begin {example} \label {ex-reltropcurve}
  The following curve contributes to $N^{\trop}_{5, (6,1,0,1)}((1),(2,1))$ with
  multiplicity $\frac{1}{2}\cdot 2\cdot 2=2$. We have drawn a grey dot at the
  end of the up most left end in order to indicate that its $y$-coordinate is
  fixed.

  \begin {center} \input {pics/expsicurve2} \end {center}
\end{example}

\begin {remark}[The equality $ \tilde N^{\trop}_{d,\kk}(\alpha,\beta) =
    \tilde N_{d,\kk}(\alpha,\beta)$] \label {rem-eqreltropclass}
  There is no direct correspondence known between the numbers
  $\tilde{N}^{\trop}_{d, {\bf k}}(\alpha, \beta)$ and $\tilde{N}_{d, {\bf
  k}}(\alpha, \beta)$. However, we prove in Theorem \ref {thm:relativepsicount}
  that $N^{\trop}_{d, {\bf k}}(\alpha, \beta)=N^{\floor}_{d,
  {\bf k}}(\alpha, \beta)$, and we show in Theorem \ref {thm-CHfloor} and
  Corollary \ref{cor-CHclass} that the numbers $N^{\floor}_{d, {\bf
  k}}(\alpha, \beta)$ and $N_{d, {\bf k}}(\alpha, \beta)$ satisfy the
  same recursive relation. It follows that $N^{\floor}_{d, {\bf
  k}}(\alpha, \beta)=N_{d, {\bf k}}(\alpha, \beta)$ and thus also that
  $N^{\trop}_{d, {\bf k}}(\alpha, \beta)=N_{d, {\bf k}}(\alpha,
  \beta)$. Of course, the analogous statements hold for the numbers
  $\tilde{N}^{\trop}_{d, {\bf k}}(\alpha, \beta)$ and $\tilde{N}_{d,
  {\bf k}}(\alpha, \beta)$ as well.
\end {remark}

%% file: pics/expsicurve.tex
\begin{picture}(0,0)%
\includegraphics{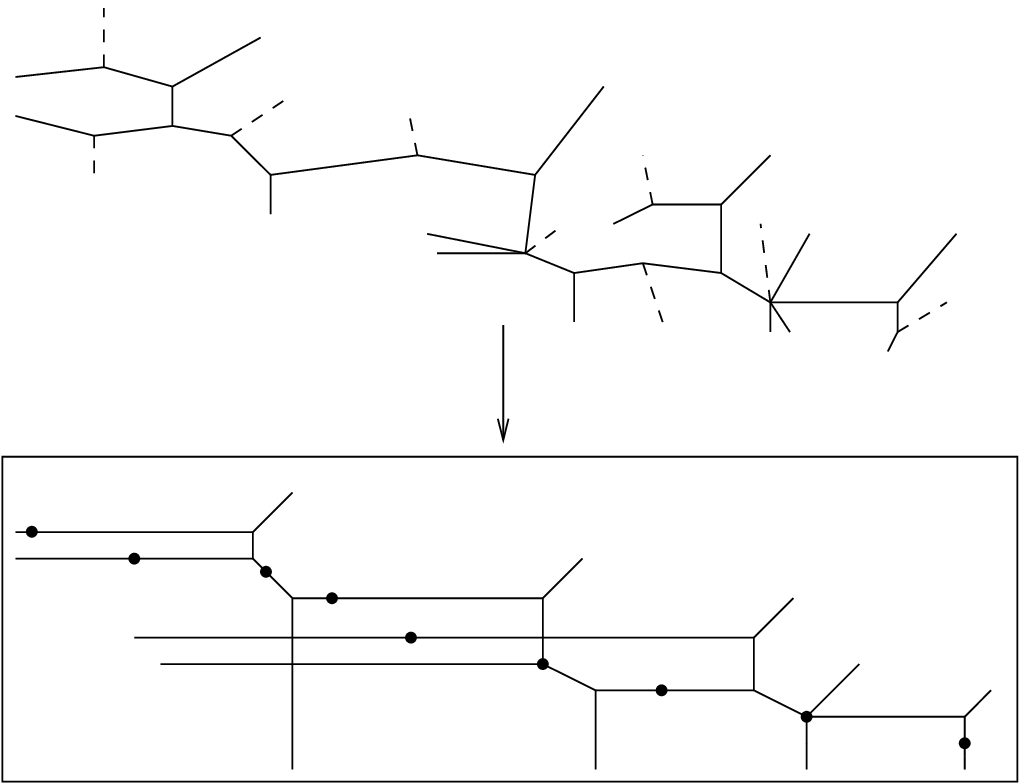}%
\end{picture}%
\setlength{\unitlength}{3947sp}%
\begingroup\makeatletter\ifx\SetFigFont\undefined%
\gdef\SetFigFont#1#2#3#4#5{%
  \reset@font\fontsize{#1}{#2pt}%
  \fontfamily{#3}\fontseries{#4}\fontshape{#5}%
  \selectfont}%
\fi\endgroup%
\begin{picture}(4896,3940)(2764,-6448)
\put(2839,-5113){\makebox(0,0)[lb]{\smash{{\SetFigFont{10}{12.0}{\familydefault}{\mddefault}{\updefault}{\color[rgb]{0,0,0}$h(x_1)$}%
}}}}
\put(6052,-5970){\makebox(0,0)[lb]{\smash{{\SetFigFont{8}{9.6}{\familydefault}{\mddefault}{\updefault}{\color[rgb]{0,0,0}$2$}%
}}}}
\put(5364,-6039){\makebox(0,0)[rb]{\smash{{\SetFigFont{10}{12.0}{\familydefault}{\mddefault}{\updefault}{\color[rgb]{0,0,0}$h(x_6)$}%
}}}}
\put(6582,-6278){\makebox(0,0)[rb]{\smash{{\SetFigFont{10}{12.0}{\familydefault}{\mddefault}{\updefault}{\color[rgb]{0,0,0}$h(x_8)$}%
}}}}
\put(7576,-5068){\makebox(0,0)[rb]{\smash{{\SetFigFont{10}{12.0}{\familydefault}{\mddefault}{\updefault}{\color[rgb]{0,0,0}$\RR^2$}%
}}}}
\put(3852,-4245){\makebox(0,0)[lb]{\smash{{\SetFigFont{10}{12.0}{\familydefault}{\mddefault}{\updefault}{\color[rgb]{0,0,0}$\Gamma$}%
}}}}
\put(5244,-4561){\makebox(0,0)[lb]{\smash{{\SetFigFont{10}{12.0}{\familydefault}{\mddefault}{\updefault}{\color[rgb]{0,0,0}$h$}%
}}}}
\put(4177,-3107){\makebox(0,0)[b]{\smash{{\SetFigFont{10}{12.0}{\familydefault}{\mddefault}{\updefault}{\color[rgb]{0,0,0}$x_3$}%
}}}}
\put(5844,-3376){\makebox(0,0)[b]{\smash{{\SetFigFont{10}{12.0}{\familydefault}{\mddefault}{\updefault}{\color[rgb]{0,0,0}$x_5$}%
}}}}
\put(5960,-4376){\makebox(0,0)[b]{\smash{{\SetFigFont{10}{12.0}{\familydefault}{\mddefault}{\updefault}{\color[rgb]{0,0,0}$x_7$}%
}}}}
\put(6366,-3697){\makebox(0,0)[lb]{\smash{{\SetFigFont{10}{12.0}{\familydefault}{\mddefault}{\updefault}{\color[rgb]{0,0,0}$x_8$}%
}}}}
\put(3266,-2655){\makebox(0,0)[b]{\smash{{\SetFigFont{10}{12.0}{\familydefault}{\mddefault}{\updefault}{\color[rgb]{0,0,0}$x_1$}%
}}}}
\put(4717,-3170){\makebox(0,0)[b]{\smash{{\SetFigFont{10}{12.0}{\familydefault}{\mddefault}{\updefault}{\color[rgb]{0,0,0}$x_4$}%
}}}}
\put(5462,-3749){\makebox(0,0)[b]{\smash{{\SetFigFont{10}{12.0}{\familydefault}{\mddefault}{\updefault}{\color[rgb]{0,0,0}$x_6$}%
}}}}
\put(7354,-4189){\makebox(0,0)[lb]{\smash{{\SetFigFont{10}{12.0}{\familydefault}{\mddefault}{\updefault}{\color[rgb]{0,0,0}$x_9$}%
}}}}
\put(3218,-3674){\makebox(0,0)[b]{\smash{{\SetFigFont{10}{12.0}{\familydefault}{\mddefault}{\updefault}{\color[rgb]{0,0,0}$x_2$}%
}}}}
\end{picture}%

%% file: pics/expsicurve1.tex
\begin{picture}(0,0)%
\includegraphics{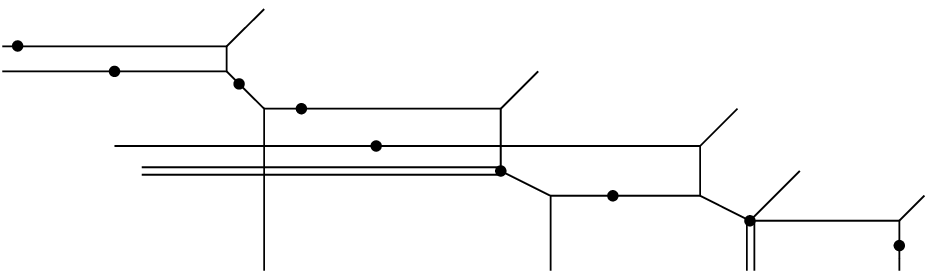}%
\end{picture}%
\setlength{\unitlength}{3947sp}%
\begingroup\makeatletter\ifx\SetFigFont\undefined%
\gdef\SetFigFont#1#2#3#4#5{%
  \reset@font\fontsize{#1}{#2pt}%
  \fontfamily{#3}\fontseries{#4}\fontshape{#5}%
  \selectfont}%
\fi\endgroup%
\begin{picture}(4449,1279)(3964,-6373)
\put(6690,-5964){\makebox(0,0)[lb]{\smash{{\SetFigFont{8}{9.6}{\familydefault}{\mddefault}{\updefault}{\color[rgb]{0,0,0}$2$}%
}}}}
\end{picture}%

%% file: pics/expsicurve2.tex
\begin{picture}(0,0)%
\includegraphics{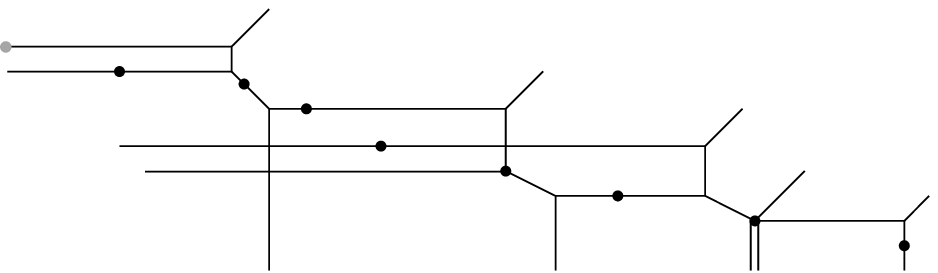}%
\end{picture}%
\setlength{\unitlength}{3947sp}%
\begingroup\makeatletter\ifx\SetFigFont\undefined%
\gdef\SetFigFont#1#2#3#4#5{%
  \reset@font\fontsize{#1}{#2pt}%
  \fontfamily{#3}\fontseries{#4}\fontshape{#5}%
  \selectfont}%
\fi\endgroup%
\begin{picture}(4471,1278)(3940,-6371)
\put(5742,-6010){\makebox(0,0)[lb]{\smash{{\SetFigFont{8}{9.6}{\familydefault}{\mddefault}{\updefault}{\color[rgb]{0,0,0}$2$}%
}}}}
\put(6712,-5964){\makebox(0,0)[lb]{\smash{{\SetFigFont{8}{9.6}{\familydefault}{\mddefault}{\updefault}{\color[rgb]{0,0,0}$2$}%
}}}}
\end{picture}%

%% file: psifloor.tex
\section{Psi-floor diagrams} \label {sec-floor}

\subsection {Absolute Psi-floor diagrams}

Floor diagrams, introduced by Brugall\'e and Mikhalkin \cite{BM1,BM2}, are
enriched directed graphs which, if counted correctly, enumerate plane curves
satisfying certain point and tangency conditions. In the following, we
generalize this definition to Psi-floor diagrams, and prove that they enumerate
tropical plane curves satisfying point, tangency, and Psi-conditions. Let us
begin with an example motivating in which sense floor diagrams extract the
combinatorial essence of a tropical curve. Return to Example \ref
{ex-nonrelativecurve}. There we have already chosen a horizontally stretched
configuration (see Definition 3.1 of \cite{FM}, they use vertically stretched).
So we expect the tropical curve to decompose into floors, and the floors are
connected by horizontal edges only. Let us point this out in the example:

\begin {center} \input {pics/floors} \end {center}

Each floor is fixed by one point, and the horizontal edges which are not
adjacent to a Psi-point are also fixed by a point. We can already see that the
presence of points satisfying Psi-conditions may lead to multiple floors ---
the second floor from the right is of degree $2$, since it contains two ends of
direction $(0,-1)$ resp.\ $(1,1)$. The marked Psi-floor diagram of this curve
can be found in step 3 of Definition \ref{def-marking}.

In the original setting of floor diagrams \cite{BM1,BM2,FM} there are only
single floors with one end of direction $(0,-1)$ and one of direction $(1,1)$.
There the idea is to shrink each floor to one vertex, and then first consider
a weighted graph on the vertex set of all floors (a floor diagram). The weights
of the edges correspond to the weights of the corresponding edges of the
tropical curve. One obtains the ``marking'' of the floor diagram by adding in
the ends and points on horizontal edges. Since any direction vector of an edge
inside a floor has $y$-coordinate $1$, a horizontal edge of weight $i$ has to
end at two vertices of multiplicity $i$ each. Therefore, the multiplicity of a
floor diagram equals the product over the squares of these weights.

Our setting is similar, but differs in a few features which we address now
before giving the precise definition. We have seen already that multiple floors
can occur. Consider a contracted end with Psi-condition $\psi^a$ in a multiple
floor of degree $d'$ (i.e.\ $d'$ ends of direction $(0,-1)$ resp.\ $(1,1)$
belong to the floor). If we remove the contracted end from the abstract graph,
we produce $a+2$ connected components. Therefore, we must have $a+2\geq 2d'$
(the \emph{string inequality}), since otherwise there would be a connected
component which contains two ends, and thus a string (see Definition 3.5 of
\cite{GM053}), in contradiction to the general position of the points.

As explained above, a multiple floor of degree $d'$ has $d'$ ends of direction
$(0,-1)$ and $(1,1)$. Furthermore, it has some ``incoming edges'' of directions
$(-m,0)$ and some ``outgoing edges'' of directions $(m,0)$ (for some $ m \in
\ZZ_{>0} $). Thus the balancing condition for the $x$-coordinate implies that
the sum of the weights of the incoming edges equals the sum of the weights of
the outgoing edges plus $d'$. This will be called the \emph {divergence
condition} of the floor diagram. Note however that we do not draw left ends of
the tropical curve in the floor diagram. Therefore the divergence condition
will be an inequality (that determines how many left ends are adjacent to a
floor) and not an equality.

Psi-points do not need to lie on floors --- they can also lie on horizontal
edges, as the following picture shows.

\begin {center} \input {pics/zerofloor} \end {center}

Since there may be bounded edges from other floors adjacent to such a Psi-point
on a horizontal edge, we have to include these points in the underlying floor
diagram. Therefore, we introduce degree-$0$ vertices corresponding to these
points. As we do not draw ends in the floor diagram, the valence of such a
degree zero vertex has to be the correct one after adding the ends. The
Psi-floor diagram (for details see below) of the tropical curve above is

\begin{center}
  \begin{picture}(20,30)(40,-15)\setlength{\unitlength}{5pt}\thicklines
  \Eeee\eEee
  \put(0,0){\circle*{2}}
  \put(10,0){\circle*{2}}
  \put(20,0){\circle*{2}}
  \put(0,-3){\makebox(0,0){$1 \, 2$}}
  \put(10,-3){\makebox(0,0){$0 \, 1$}}
  \put(20,-3){\makebox(0,0){$2 \, 2$}}
  \put(15,1.5){\makebox(0,0){$2$}} 
  \put(5,0){\vector(-1,0){1}} 
  \put(15,0){\vector(-1,0){1}} 
  \put(22,0){\makebox(0,0){$.$}}
  \end{picture}
\end{center}

Here is the formal definition:

\begin{definition} \label{def:floordiagram}
  A \emph{(rational) Psi-floor diagram} $\D$ is a connected, directed graph
  $(V,E)$ of genus $0$ on a linearly ordered vertex set $(V,<)$ with edge
  weights $ \omega(e) \in \ZZ_{>0} $ for all edges $e \in E$, together with
  pairs $ (d_v, a_v) \in \ZZ_{\ge 0}^2$ for each vertex $v$ in $V$ (which we
  call the \emph{degree} $d_v$ and the \emph{Psi-power} $a_v$ of $v$),
  satisfying:
  \begin{enumerate}
  \item The edge directions preserve the vertex order, i.e.\ for every edge $
    v \to w$ we have $v < w$.
  \item There are no edges between degree-$0$ vertices, i.e.\ if $v \to w$ is
    an edge then $d_v >0$ or $d_w > 0$.
  \item For each $v \in V$ at least one of the numbers $d_v$ and $a_v$ is
    positive.
  \item For each $v \in V$ we have $a_v - 2(d_v -1) \ge 0$ (string inequality).
  \item (Divergence condition) For every vertex $v$ we have
    \begin{displaymath}
      \dive(v)
        := \sum_{\substack {
             \text {edges $e$} \\[0.5ex]
             v \, \stackrel e\to \, w}
         } \omega(e) -
           \sum_{\substack {
             \text {edges $e$} \\[0.5ex]
             w \, \stackrel e\to \, v}
           } \omega(e)
	\leq d_v.
    \end{displaymath}
    This means that at every vertex of $\D$ the total weight of the outgoing
    edges is larger by at most $d_v$ than the total weight of the incoming
    edges.
  \item \label{itm:valence}
    If $d_v = 0$ for a vertex $v$, then $\val(v) = a_v + 2 + \dive(v)$ (where
    $\val(v)$ is the valence of $v$).  
  \end{enumerate}
\end{definition}

We call $d(\D) = \sum_{v \in V} d_v$ the \emph{degree} of a Psi-floor diagram
$\D$. A \emph{floor} of $\D$ is a vertex of positive degree. The \emph{type} of
$\D$ is ${\bf k}(\D) = ({\bf k}_0, {\bf k}_1, \dots)$, where ${\bf k}_i$ is the
number of vertices $v$ of $\D$ with $a_v = i$ for all $i \ge 1$, and $\kk_0$ is
the number of vertices $v$ with $a_v = 0$ plus $3d-1-I\kk-\#V$. The number
$3d-1-I\kk-\#V$ that we add to $\kk_0$ equals the number of vertices of
Psi-power $0$ that we will add later and which makes the equality
$I\kk=3d-1-|\kk|$ hold. The \emph{multiplicity} $\mu(\D)$ of $\D$ is given by
\begin{displaymath}
  \mu(\D) := \prod_{\text{edges }e} \omega(e)^2 \prod_{\substack {
    v  \, \stackrel e\to \, w \\[0.5ex]
    \text{s.t.\ } d_v = 0 \\[0.5ex]
    \text{or } d_w = 0
  }} \frac{1}{\omega(e)} \prod_{v: \, d_v = 0} \frac{1}{|\dive(v)| \, !}.
\end{displaymath}
The first factor in the definition of multiplicity corresponds, as in the
original definition of floor diagram, to vertices adjacent to edges of higher
weight. If an edge of higher weight is adjacent to a contracted end however
(e.g.\ at a vertex of degree $0$), this vertex does not contribute and so we
have to divide out by one factor of $\omega(e)$ again. The last factor
contributes to the factor $\nu_C$ in Definition \ref {def-mult} of the
multiplicity of a tropical curve, which arises because ends of the same
direction are adjacent to a vertex.

We draw Psi-floor diagrams using the convention that vertices in increasing
order are arranged left to right, thereby adopting the convention of \cite{FM}.
Note that in this paper we draw the corresponding tropical curves in the
\emph{opposite} direction. We write the pair $(d_v, a_v)$ below each vertex
$v$. Edge weights of $1$ are omitted.

\begin{example} \label{ex:floordiagram}
  An example of a Psi-floor diagram $\D$ of degree $d = 5$, type ${\bf k} =
  (7,0,1,1)$, divergences $1,1,-1,-1$, and multiplicity $\mu(\D) = 4$ is drawn
  below.

\begin{displaymath}
    \begin{picture}(50,35)(48,-15)\setlength{\unitlength}{5pt}\thicklines
    \oooo\Eeee\eEee\eeEe
    \put(15,1.5){\makebox(0,0){$2$}} 
    \put(5,0){\vector(1,0){1}} 
    \put(15,0){\vector(1,0){1}} 
    \put(25,0){\vector(1,0){1}} 
    \put(0,-3){\makebox(0,0){$1 \, 0$}}
    \put(10,-3){\makebox(0,0){$2 \, 3$}}
    \put(20,-3){\makebox(0,0){$1 \, 2$}}
    \put(30,-3){\makebox(0,0){$1 \, 0$}}
    \end{picture}
  \end{displaymath}
\end{example}

Given a Psi-floor diagram $\D$ we define, for every floor $v$, the sets $I(v)$
and $O(v)$ by
\begin{displaymath} \begin{split}
  I(v) & := \{w \to v : d_w > 0 \}, \\
  O(v) & := \{v \to w: d_w > 0 \} \cup \coprod
    \{ v \stackrel{1}{\to} \circ \},
\end{split} \end{displaymath}
where the latter set is a disjoint union of the outgoing edges of $\D$ at $v$
augmented by $d_v - \dive(v)$ many \emph{indistinguishable} edges of weight $1$
directed away from $v$ ending in distinct vertices $\circ$. These
indistinguishable extra ends correspond to left ends of the tropical curve
starting at this floor.

{\bf Example \ref{ex:floordiagram} (continued).}
  We draw the sets $I(v)$ and $O(v)$ by augmenting the Psi-floor diagrams at
  the respective vertices. If, for example, $v$ is the third black vertex from
  the left, then $O(v)$ consists of the edge between $v$ and the fourth black
  vertex and the two edges of weight $1$ connecting $v$ with the two adjacent
  white vertices.

  \begin{displaymath}
    \begin{picture}(50,45)(65,-20)\setlength{\unitlength}{5pt}\thicklines
    \oooo\Eeee\eEee\eeEe
    \put(15,1.5){\makebox(0,0){$2$}} 
    \put(5,0){\vector(1,0){1}} 
    \put(15,0){\vector(1,0){1}} 
    \put(25,0){\vector(1,0){1}} 
    \put(0,-3){\makebox(0,0){$1 \, 0$}}
    \put(10,-3){\makebox(0,0){$2 \, 3$}}
    \put(20,-3){\makebox(0,0){$1 \, 2$}}
    \put(30,-3){\makebox(0,0){$1 \, 0$}}
    \put(10,+0.0){\line(1,1){4.2}}
    \put(10.7,+0.7){\vector(1,+1){3}}
    \put(15,+5){\circle{2}}
    \put(20,+0.0){\line(1,1){4.2}}
    \put(20.7,+0.7){\vector(1,+1){3}}
    \put(25,+5){\circle{2}}
    \put(20,+0){\line(2,1){8.4}}
    \put(20,+0){\vector(2,+1){6}}
    \put(29,+5){\circle{2}}
    \put(30,+0.0){\line(1,1){4.2}}
    \put(30.7,+0.7){\vector(1,+1){3}}
    \put(35,+5){\circle{2}}
    \put(30,+0.0){\line(2,1){8.4}}
    \put(30,+0.0){\vector(2,+1){6}}
    \put(39,+5){\circle{2}}
    \end{picture}
  \end{displaymath}

An \emph {edge choice} is a collection $\C(\D)$ of subsets $C(v) \subset I(v)
\cup O(v)$, one for each floor $v$ of $\D$, satisfying $|C(v)| = a_v -2 (d_v
-1)$, and such that $C(v) \cap C(w) = \emptyset$ for distinct floors $v$ and
$w$. If $d_v = 0$ for a vertex $v$ we set $C(v) = \emptyset$. The
\emph{local multiplicity at} $v$ of such a choice is
\begin{displaymath}
  \mu_{v, C(v)} := \begin {cases}
    \frac {d_v^{i(v)}}{d_v!} \cdot \frac{d_v^{o(v)}}{d_v!}
      & \mbox {if $ d_v > 0 $}, \\[0.5ex]
    1 & \mbox {if $ d_v = 0 $}.
  \end {cases}
\end{displaymath}
where $i(v)  = |I(v) \backslash C(v)|$ and $o(v)  = |O(v) \backslash C(v)|$ are
the number of non-chosen edges in $I(v)$ and $O(v)$, respectively.

The chosen edges will later correspond to the edges of the tropical curve that
are directly adjacent to the Psi-point; the non-chosen edges to those belonging
to the floor but not directly adjacent to the Psi-point. We will see later in
Lemma \ref {lem-stirling} and the proof of Theorem \ref {thm:psicount} that the
local multiplicity at $v$ of an edge choice takes the possibilities for the
degree-$d_v$ floor and the contribution to the multiplicity $\nu_C$ of
Definition \ref {def-mult} into account. 

The \emph{multiplicity} $\mu(\C)$ of the edge choice $\C(\D)$ of the Psi-floor
diagram $\D$ is
\begin{displaymath}
  \mu(\C)
    := \prod_{v \in V} \mu_{v,C(v)} \, \frac 1{|C(v) \cap \{v \to \circ \}|\,!}
       \, \prod_{e \in C(v)} \frac{1}{\omega(e)}.
\end{displaymath}
As before, the multiplicity of an edge choice takes for each floor a
combination of contributions to $\nu_C$ and possibilities for a floor into
account, furthermore additional contributions to $\nu_C$ and factors of
$\frac{1}{\omega(e)}$ that arise because an edge of weight $\omega(e)$ is
adjacent to a contracted end.  

\begin{example} \label{ex:edgechoice}
  We picture an edge choice $\C(\D)$ by thickening all edges in $C(v)$ at $v$,
  for all vertices $v$ of $\D$. Below is an edge choice for the Psi-floor
  diagram of Example \ref{ex:floordiagram}. Its multiplicity is $\mu(\C) =
  \tfrac 12 $. Notice that $|C(v)| = a_v - 2(d_v-1) $ for all $v$ since none of
  the vertices has degree zero. 

  \begin{displaymath}
    \begin{picture}(50,40)(65,-15)\setlength{\unitlength}{5pt}\thicklines
    \oooo\Eeee\eEee\eeEe
    \put(15,1.5){\makebox(0,0){$2$}} 
    \put(5,0){\vector(1,0){1}} 
    \put(15,0){\vector(1,0){1}} 
    \put(25,0){\vector(1,0){1}} 
    \put(0,-3){\makebox(0,0){$1 \, 0$}}
    \put(10,-3){\makebox(0,0){$2 \, 3$}}
    \put(20,-3){\makebox(0,0){$1 \, 2$}}
    \put(30,-3){\makebox(0,0){$1 \, 0$}}
    \put(10,+0.0){\line(1,1){4.2}}
    \put(10.7,+0.7){\vector(1,+1){3}}
    \put(15,+5){\circle{2}}
    \put(20,+0.0){\line(1,1){4.2}}
    \put(20.7,+0.7){\vector(1,+1){3}}
    \put(25,+5){\circle{2}}
    \put(20,+0){\line(2,1){8.4}}
    \put(20,+0){\vector(2,+1){6}}
    \put(29,+5){\circle{2}}
    \put(30,+0){\line(1,1){4.2}}
    \put(30.7,+0.7){\vector(1,+1){3}}
    \put(35,+5){\circle{2}}
    \put(30,+0){\line(2,1){8.4}}
    \put(30,+0){\vector(2,+1){6}}
    \put(39,+5){\circle{2}}
    \linethickness{1mm}
    \put(6,0){\line(1,0){4}}
    \put(20,0){\qbezier(0,0)(1,1)(2.5,2.5)}
    \put(20,0){\qbezier(0,0)(2,1)(4,2)}
    \end{picture}
  \end{displaymath}
\end{example}

\begin{definition} \label{def-marking}
  A \emph{marking} of a Psi-floor diagram $\D$ with an edge choice $\C$ is
  defined by the following three-step process which we will illustrate in the
  case of Example \ref{ex:edgechoice}.

  {\bf Step 1:} For each vertex $v$ of $\D$ create $d_v- \dive(v) - |C(v)\cap
  \{v \to \circ \}|$ many new vertices in $\D$ and connect them to $v$ with new
  edges directed away from $v$.

  \begin{displaymath}
    \begin{picture}(50,45)(65,-18)\setlength{\unitlength}{5pt}\thicklines
    \oooo\Eeee\eEee\eeEe
    \put(15,1.5){\makebox(0,0){$2$}} 
    \put(5,0){\vector(1,0){1}} 
    \put(15,0){\vector(1,0){1}} 
    \put(25,0){\vector(1,0){1}} 
    \put(0,-3){\makebox(0,0){$1 \, 0$}}
    \put(10,-3){\makebox(0,0){$2 \, 3$}}
    \put(20,-3){\makebox(0,0){$1 \, 2$}}
    \put(30,-3){\makebox(0,0){$1 \, 0$}}
    \put(10,+0.0){\line(1,1){4.2}}
    \put(10.7,+0.7){\vector(1,+1){3}}
    \put(15,+5){\circle{2}}
    \put(30,+0){\line(1,1){4.2}}
    \put(30.7,+0.7){\vector(1,+1){3}}
    \put(35,+5){\circle{2}}
    \put(30,+0){\line(2,1){8.4}}
    \put(30,+0){\vector(2,+1){6}}
    \put(39,+5){\circle{2}}
    \linethickness{1mm}
    \put(6,0){\line(1,0){4}}
    \end{picture}
  \end{displaymath}

  These correspond exactly to the non-chosen edges $ v \to \circ $ above,
  i.e.\ to the left ends of the tropical curve that are not directly adjacent
  to the Psi-point in the floor (and therefore have to be fixed later by a
  point condition).

  {\bf Step 2:} Subdivide each non-chosen edge of the original Psi-floor
  diagram $\D$ between floors into two directed edges by introducing a new
  vertex for each such edge. The new edges inherit their weights and
  orientations. Call the resulting graph $\tilde{\D}$.

  \begin{displaymath}
    \begin{picture}(50,45)(65,-18)\setlength{\unitlength}{5pt}\thicklines
    \oooo\Eeee
    \put(5,0){\vector(1,0){1}} 
    \put(10,+0.0){\line(1,0){4}}
    \put(15,+0){\circle{2}}
    \put(16,+0.0){\line(1,0){4}}
    \put(20,+0.0){\line(1,0){4}}
    \put(13,1.5){\makebox(0,0){$2$}}
    \put(18,1.5){\makebox(0,0){$2$}}
    \put(25,+0){\circle{2}}
    \put(26,+0.0){\line(1,0){4}}
    \put(12,0){\vector(1,0){1}} 
    \put(17,0){\vector(1,0){1}} 
    \put(22,0){\vector(1,0){1}} 
    \put(27,0){\vector(1,0){1}} 
    \put(0,-3){\makebox(0,0){$1 \, 0$}}
    \put(10,-3){\makebox(0,0){$2 \, 3$}}
    \put(20,-3){\makebox(0,0){$1 \, 2$}}
    \put(30,-3){\makebox(0,0){$1 \, 0$}}
    \put(10,+0.0){\line(1,1){4.2}}
    \put(10.7,+0.7){\vector(1,+1){3}}
    \put(15,+5){\circle{2}}
    \put(30,+0){\line(1,1){4.2}}
    \put(30.7,+0.7){\vector(1,+1){3}}
    \put(35,+5){\circle{2}}
    \put(30,+0){\line(2,1){8.4}}
    \put(30,+0){\vector(2,+1){6}}
    \put(39,+5){\circle{2}}
    \linethickness{1mm}
    \put(6,0){\line(1,0){4}}
    \end{picture}
  \end{displaymath}

  These extra vertices correspond to points on horizontal bounded edges with no
  Psi-condition.

  {\bf Step 3:} Order the vertices of $\tilde{\D}$ linearly, extending the
  order of the vertices of the original Psi-floor diagram $\D$, such that (as
  in $\D$) each edge is directed from a smaller vertex to a larger vertex.

  \begin{displaymath}
    \begin{picture}(50,40)(70,-20)\setlength{\unitlength}{5pt}\thicklines
    \put(0,+0){\circle*{2}}
    \put(5,+0){\circle*{2}}
    \put(10,+0){\circle{2}}
    \put(15,+0){\circle*{2}}
    \put(20,+0){\circle{2}}
    \put(25,+0){\circle{2}}
    \put(30,+0){\circle*{2}}
    \put(35,+0){\circle{2}}
    \put(40,+0){\circle{2}}
    \put(0,-3){\makebox(0,0){$1 \, 0$}}
    \put(5,-3){\makebox(0,0){$2 \, 3$}}
    \put(15,-3){\makebox(0,0){$1 \, 2$}}
    \put(30,-3){\makebox(0,0){$1 \, 0$}}
    \put(0,+0.0){\line(1,0){5}}
    \put(1.5,0){\vector(1,0){1}}
    \put(5,+0.0){\line(1,0){4}}
    \put(8,0){\vector(1,0){1}}
    \put(8,1.5){\makebox(0,0){$2$}} 
    \put(11,+0.0){\line(1,0){4}}
    \put(13,0){\vector(1,0){1}}
    \put(12.5,1.5){\makebox(0,0){$2$}} 
    \put(26,+0.0){\line(1,0){4}}
    \put(28,0){\vector(1,0){1}}
    \put(30,+0.0){\line(1,0){4}}
    \put(33,0){\vector(1,0){1}}
    \qbezier(5.8,0.6)(8,4)(12.5,4)\qbezier(12.5,4)(17,4)(19.2,0.6)
    \put(12.5,4){\vector(1,0){1}}
    \qbezier(15.8,-0.6)(17,-3)(20,-3)\qbezier(20,-3)(23,-3)(24.2,-0.6)
    \put(20,-3){\vector(1,0){1}}
    \qbezier(30.8,0.6)(32,3)(35,3)\qbezier(35,3)(38,3)(39.2,0.6)
    \put(35,3){\vector(1,0){1}}
    \linethickness{1mm}
    \put(3,0){\line(1,0){2}}
    \end{picture}
  \end{displaymath}

  The extended graph $\tilde{\D}$ together with the linear order on its
  vertices is called a \emph{marked Psi-floor diagram}, or a \emph{marking} of
  the Psi-floor diagram $\D$.
\end{definition}

We added $d_v-\dive(v)$ white end vertices for each $v\in V$ before picking the
edge choice. It follows by induction that altogether we add $d$ white end
vertices. However, in step 1 of Definition \ref{def-marking} we really only
add the non-chosen ones. In step 2 we subdivide each of the non-chosen edges.
There are $\#V-1$ edges, since the Psi-floor diagram is a rational graph. Thus,
altogether we add $d+\#V-1$ minus the number of chosen edges white vertices,
i.e.\  $d+\#V-1-\sum_{v\in V} (a_v-2(d_v-1)) = 3d-1-I\kk-\#V$. It follows that
$\kk_0$ equals the number of vertices $v$ of the floor diagram with $a_v=0$
plus the number of white vertices in the marking.

We want to count marked Psi-floor diagrams up to equivalence. Two such
$\tilde{\D}_1$, $\tilde{\D}_2$ are \emph{equivalent} if $\tilde{\D}_1$ can be
obtained from $\tilde{\D}_2$ by permuting edges without changing their weights,
i.e.\ if there exists an automorphism of weighted graphs which preserves the
vertices of $\D$ and maps $\tilde{\D}_1$ to $\tilde{\D}_2$. 

The \emph{number of markings} $\nu(\D,\C)$ is the number of marked Psi-floor
diagrams $\tilde{\D}$ up to equivalence. In the example, we have $ \nu(\D,\C)=7
$: the white $1$-valent vertex adjacent to the second black vertex (counted
from the left) can be inserted in $2$ ways between the second and third black
vertex, in $2$ ways between the third and fourth black vertex, and in $3$ ways
right of the fourth black vertex.

By specializing to the case $a_v = 0$ for all vertices $v$ of $\D$ we recover
the definition of labeled floor diagrams and their markings of Fomin and
Mikhalkin \cite{FM}. In this case all floors necessarily have degree $d_v = 1$
and no edges get chosen (so $C(v) = \emptyset$ for all vertices $v$).

\begin {definition}[$ N^{\floor}_{d,\kk} $ and $ \tilde N^{\floor}_{d,\kk} $]
  Let $d \ge 1$ and ${\bf k}$ be a sequence of non-negative integers with $I
  {\bf k} = 3d - 1 - |{\bf k}|$. Set
  \begin{displaymath}
    N^{\floor}_{d, {\bf k}}
      := \sum_\D \mu(\D) \,
         \sum_\C \mu(\C) \, \nu(\D, \C), 
  \end{displaymath}
  where the first sum is over all Psi-floor diagrams of degree $d$ and type
  ${\bf k}$, and the second sum is over all edge choices $\C$ of $\D$.
  Correspondingly (see Definition \ref {def-ndk}), we set $ \tilde
  N^{\floor}_{d,\kk} := \frac {\kk!}{|\kk|\,!} \, N^{\floor}_{d,\kk} $.
\end{definition}

\begin{remark} \label{rem-nfloortilde}
  We can also define the numbers $\tilde{N}^{\floor}_{d, {\bf k}}$ directly
  using Psi-floor diagrams by requiring that the Psi-powers of the vertices of
  the marked Psi-floor diagram (the Psi-powers of the white vertices that are
  not present in the underlying Psi-floor diagram have Psi-power $0$) occur in
  a particular order, and by marking the white end vertices with numbers from
  $1$ to $d$.
\end{remark}

\begin{example}
  As an example in degree $d = 4$ we compute the number
  \begin{displaymath}
    \tilde{N}^{\floor}_{4,(1,0,0,0,2)} = \frac{1}{4}.
  \end{displaymath}
  There are three markings of Psi-floor diagrams of degree $4$ and type
  $(1,0,0,0,2)$ which have the Psi-powers in the order $(a_1,a_2,a_3) =
  (0,4,4)$. (Remember that we draw Psi-floor diagrams from left to right and
  therefore need to invert the order of the $a_i$.) Every other order of the
  $a_i$ yields the same answer.

  \begin{displaymath}
    \begin{picture}(0,80)(200,-50)\setlength{\unitlength}{2.5pt}\thicklines
    \put(0,0){\circle*{4}}
    \put(20,0){\circle*{4}}
    \put(40,0){\circle{4}}
    \put(0,-6){\makebox(0,0){$2 \, 4$}}
    \put(20,-6){\makebox(0,0){$2 \, 4$}}
    \put(2,0){\line(1,0){16}}
    \put(2,0){\vector(1,0){13}}
    \put(22,0){\line(1,0){16}}
    \put(22,0){\vector(1,0){13}}
    \put(10,10){\circle{4}}
    \put(0,+0.0){\line(1,1){8.4}}
    \put(1.4,+1.4){\vector(1,+1){6}}
    \put(20,+0){\line(1,1){8.4}}
    \put(21.4,+1.4){\vector(1,+1){6}}
    \put(30,+10){\circle{4}}
    \put(20,+0){\line(2,1){16.8}}
    \put(20,+0){\vector(2,+1){14.4}}
    \put(38,+10){\circle{4}}
    \linethickness{1.5mm}
    \put(1.5,0){\line(1,0){8}}
    \put(0,0){\qbezier(0,0)(2,2)(5,5)}
    \put(20,0){\qbezier(0,0)(2,2)(5,5)}
    \put(20,0){\qbezier(0,0)(4,2)(8,4)}
    \put(20,-14){\makebox(0,0){
      $\frac{1}{2} \cdot \frac{1}{2} \cdot \frac{1}{2!} = \frac{1}{8}$}}
    \end{picture}
    \begin{picture}(0,80)(50,-50)\setlength{\unitlength}{2.5pt}\thicklines
    \put(0,0){\circle*{4}}
    \put(20,0){\circle*{4}}
    \put(40,0){\circle{4}}
    \put(0,-6){\makebox(0,0){$3 \, 4$}}
    \put(20,-6){\makebox(0,0){$1 \, 4$}}
    \put(2,0){\line(1,0){16}}
    \put(2,0){\vector(1,0){8}}
    \put(10,4){\makebox(0,0){$3$}}
    \put(22,0){\line(1,0){16}}
    \put(22,0){\vector(1,0){8}}
    \put(20,0){\line(0,1){11}}
    \put(20,2){\vector(0,1){8}}
    \put(20,13){\circle{4}}
    \put(21.4,1.4){\line(1,1){10.2}}
    \put(21.4,1.4){\vector(1,+1){8.2}}
    \put(33,13){\circle{4}}
    \put(21.4,1.4){\line(2,1){17.8}}
    \put(21.4,1.4){\vector(2,+1){14.4}}
    \put(41,11){\circle{4}}
    \linethickness{1.5mm}
    \put(18.5,0){\line(-1,0){7}}
    \put(20,0){\qbezier(0,0)(0,2)(0,6)}
    \put(21.4,1.4){\qbezier(0,0)(2,2)(5,5)}
    \put(21.4,1.4){\qbezier(0,0)(4,2)(8,4)}
    \put(20,-14){\makebox(0,0){$9 \cdot \frac{1}{3} \cdot \frac{3}{3!} \cdot
      \frac{3}{3!} \cdot \frac{1}{3!} = \frac{1}{24}$}}
    \end{picture}
    \begin{picture}(0,80)(-100,-50)\setlength{\unitlength}{2.5pt}\thicklines
    \put(0,0){\circle*{4}}
    \put(20,0){\circle*{4}}
    \put(40,0){\circle{4}}
    \put(0,-6){\makebox(0,0){$3 \, 4$}}
    \put(20,-6){\makebox(0,0){$1 \, 4$}}
    \put(2,0){\line(1,0){16}}
    \put(2,0){\vector(1,0){8}}
    \put(10,4){\makebox(0,0){$2$}}
    \put(0,0){\qbezier(0,0)(6,-9)(20,-9) \qbezier(20,-9)(32,-9)(38.4,-1.6)}
    \put(20,-9){\vector(1,0){1}}
    \put(20,0){\line(0,1){11}}
    \put(20,2){\vector(0,1){8}}
    \put(20,13){\circle{4}}
    \put(21.4,1.4){\line(1,1){10.2}}
    \put(21.4,1.4){\vector(1,+1){8.2}}
    \put(33,13){\circle{4}}
    \put(21.4,1.4){\line(2,1){17.8}}
    \put(21.4,1.4){\vector(2,+1){14.4}}
    \put(41,11){\circle{4}}
    \linethickness{1.5mm}
    \put(18.5,0){\line(-1,0){7}}
    \put(20,0){\qbezier(0,0)(0,2)(0,6)}
    \put(21.4,1.4){\qbezier(0,0)(2,2)(5,5)}
    \put(21.4,1.4){\qbezier(0,0)(4,2)(8,4)}
    \put(20,-14){\makebox(0,0){$4 \cdot \frac{1}{2} \cdot\frac{1}{3!} \cdot
      \frac{9}{3!} \cdot \frac{1}{3!} = \frac{1}{12}$}}
    \end{picture}
  \end{displaymath}

  The contribution of the third diagram, for example, arises as follows: The
  underlying Psi-floor diagram has multiplicity $4$. Choosing the weight-$2$
  edge amounts to a factor of $\tfrac{1}{2}$. The degree-$3$ vertex has no
  non-chosen incoming edges and $2$ non-chosen outgoing edges, hence the local
  multiplicity at this vertex is given by the next two factors. Lastly,
  as all $3$ additional edges at the second vertex are chosen, we need to
  multiply by $\tfrac{1}{3!}$.
\end{example}

\subsection{The equality $N^{\floor}_{d, {\bf k}}= N^{\trop}_{d, {\bf k}}$}

\begin{theorem} \label{thm:psicount}
  Let $d \ge 1$ and ${\bf k}$ be a sequence of non-negative integers with $I
  {\bf k} = 3d - 1 - |{\bf k}|$. Then $N^{\floor}_{d, {\bf k}}= N^{\trop}_{d,
  {\bf k}} $.
\end{theorem}

For the proof of Theorem \ref{thm:psicount} we need the following lemma. For
positive integers $a$ and $b$, let $S(a,b)$ denote the \emph{Stirling number of
the second kind}, i.e.\ the number of ways of partitioning an $a$-element set
into $b$ non-empty parts.

\begin{lemma}[\cite{St}, (24d)]\label{lem-stirling}
  For integers $ e, f \ge 0 $ it holds that
  \begin{displaymath}
    \sum_{0 \le g \le f} \frac{S(e,g)}{(f-g)!} = \frac{f^e}{f!}.
  \end{displaymath}
\end{lemma}

\begin{proof}[Proof of Theorem \ref{thm:psicount}]
  Pick a horizontally stretched configuration of $|\kk|$ points (see Definition
  3.1 of \cite{FM}). Our strategy is as follows: let $T$ be the set of tropical
  curves of degree $d$ satisfying the conditions, and let $F$ be the set of
  marked floor diagrams of degree $d$ and type $ \kk $. We will define a
  (surjective) map from $T$ to $F$. Let $r$ be the number of inverse images
  of a given marked floor diagram $\tilde{\D}$ in $F$. We will show that each
  such inverse image is a tropical curve $C$ of the same multiplicity $
  \mult(C)$, and that $\mult(C)\cdot r = \mu(\D)\cdot \mu(\C)$, where $\D$
  denotes the underlying floor diagram for $\tilde{\D}$ and $\C$ denotes its
  choice of edges. Of course, this will then prove the lemma.

  Consider a tropical curve in $T$; we will now explain how to construct the
  corresponding marked floor diagram in $F$. As in Theorem 4.3 of \cite{GM052}
  resp.\ Section 5 of \cite{BM2} it follows that the tropical curve decomposes
  into floors in the sense that each connected component of $\Gamma$ minus the
  horizontal edges (i.e.\ each floor) is fixed by exactly one point. (A floor
  can have higher degree here.) For each floor $v$ let $d_v$ denote its number
  of ends of direction $(0,-1)$ and $a_v$ the power of Psi of the contracted
  end (i.e.\ the valence of the adjacent vertex minus 3). Shrink each floor to
  a vertex labeled with $(d_v,a_v)$. If there is a contracted end with a
  Psi-condition on a horizontal edge, also keep this as a vertex and set
  $d_v=0$, and $a_v$ the power of Psi. Let the edges of the floor diagram be
  given by the horizontal bounded edges of the tropical curve connecting the
  floors. We orient the edges towards the left ends of the curve, and reverse
  the picture (so the left ends are on the right, and edges are oriented to the
  right). Because of the general position of the points there cannot be two
  contracted ends mapped to a horizontal line --- thus there cannot be any
  edges between vertices of degree $0$. If $d_v=0$ for a vertex we know that
  the corresponding contracted end has a Psi-condition, so then $a_v>0$. Of
  course, if $a_v=0$ then we must have $d_v>0$. If there are horizontal ends
  adjacent to a contracted end on a floor resp.\ to a contracted end with
  higher Psi-condition on a horizontal edge, drop them. The other horizontal
  ends must be adjacent to a contracted end without a Psi-condition; keep the
  contracted end as a white end vertex. Also draw white vertices on horizontal
  edges for contracted ends without a Psi-condition on horizontal edges.
  Thicken the horizontal edges which are directly adjacent to a contracted end
  on a floor. A vertex of degree $0$ in the floor diagram comes from a
  contracted end with a Psi-condition, say of power $a_v$, on a horizontal
  edge. Since the tropical curve is balanced, the sum of the weights of the
  incoming horizontal edges must equal the sum of the weights of the outgoing.
  The divergence condition for degree-$0$ vertices follows. The valence must be
  $a_v+2$ (without counting the contracted end itself). We have dropped the
  ends adjacent to this vertex however, so we have $\val(v)-\dive(v)=a_v+2$.
  Now let $v$ be a vertex of the floor diagram with $d_v>0$. This vertex comes
  from a floor of the tropical curve which contains a contracted end with
  Psi-power $a_v$. If we remove the contracted end from $\Gamma$ we produce
  $a_v+2$ connected components. The floor contains $2d_v$ ends of direction
  $(0,-1)$ resp.\ $(1,1)$. These ends must belong to different connected
  components since otherwise there would be a string (see Definition 3.5 of
  \cite{GM053}) in contradiction to the general position of the points. It
  follows that $a_v+2\geq 2d_v$ (string inequality), and that $a_v+2-2d_v$
  horizontal edges are directly adjacent to the contracted end, and thus get
  chosen (including ends, which we drop). For a vertex of the floor diagram
  with $d_v>0$, the balancing condition in the $x$-coordinate tells us that the
  divergence condition holds. It follows that we have produced a marked
  Psi-floor diagram in $F$ for the curve in $T$.

  Conversely, let now $\tilde{\D}$ be a marked floor diagram in $F$; we will
  construct its inverse images in $T$. For each white vertex and for each
  vertex of degree $0$ draw horizontal edges of the appropriate weight through
  the corresponding point $p_i$. For a vertex of degree $d_v>0$ there are
  several possibilities how it can be completed to a floor of a tropical curve.
  We have seen already that --- locally around such a floor of a tropical curve
  --- removing the contracted end produces $a_v+2$ connected components of
  which $a_v+2-2d_v$ are horizontal edges and $2d_v$ are connected components
  containing one of the $2d_v$ ends of direction $(0,-1)$ resp.\ $(1,1)$. There
  are $o(v)$ non-chosen outgoing horizontal edges connected to this floor.
  Their $y$-coordinates are fixed by other conditions. Thus they are
  distinguishable in the tropical curve, even if they are of the same weight.
  These edges must belong to the connected components containing the ends of
  direction $(1,1)$. Assume that $g$ of the $d_v$ connected components
  containing the ends $(1,1)$ also contain horizontal edges, whereas $d_v-g$
  ends of direction $(1,1)$ are directly adjacent to the contracted end. Thus
  we need to partition the set of $o(v)$ horizontal non-chosen edges into $g$
  non-empty parts, corresponding to the $g$ connected components. For each such
  choice there is exactly one possibility to complete the picture to the upper
  part of a floor of a tropical curve since the $y$-coordinates of the
  horizontal edges are fixed by other points. This part of the tropical curve
  contributes a factor of $\frac{1}{(d_v-g)!}$ to the factor $\nu_C$ of the
  multiplicity of the tropical curve because of the $d_v-g$ ends of direction
  $(1,1)$ which are directly adjacent to the contracted end. Thus we can sum up
  the possibilities with their contribution to $\nu_C$ as $\frac{S(o(v),g)}{
  (d_v-g)!}$ for each $g$. Summing over all $g$, we get $\frac{d_v^{o(v)}}{
  d_v!}$ by Lemma \ref{lem-stirling}. This situation is illustrated in Example
  \ref {ex-stirling}.

  The analogous statement holds for the lower part of the floor of the tropical
  curve and the incoming horizontal edges. For any choice of $g$ and a
  partition (both for the upper and the lower part of each floor) we can
  complete the picture uniquely to a tropical curve.

  The multiplicity of the tropical curve is a product of factors contributing
  to $\nu_C$ and vertex multiplicities. We have taken care of the factors
  contributing to $\nu_C$ inside each floor already. There can still be left
  ends adjacent to the same vertex that contribute to $\nu_C$. This happens
  either if left ends are adjacent to vertices of degree $0$ in the floor
  diagram, or if they are directly adjacent to a contracted end inside a floor,
  i.e.\ chosen. For the first situation, we get a factor of
  $\frac{1}{\dive(v)!}$, for the second situation we get a factor of
  $\frac{1}{|C(v) \cap \{v \to \circ \}|\,!}$. Now let us consider the vertex
  multiplicities. We have seen already that each floor consists of components
  with one end of direction $(1,1)$ resp.\ $(0,-1)$, and horizontal edges. The
  $y$-coordinate of any direction of an edge of such a component is therefore
  $1$, and thus any vertex adjacent to a horizontal edge of weight $\omega(e)$
  is of multiplicity $\omega(e)$. If a horizontal edge is adjacent to a
  contracted end however, this vertex does not contribute. If this contracted
  end comes from a white vertex however, there is another horizontal edge of
  the same weight adjacent to it. Thus, any horizontal edge in the floor
  diagram (without the marking) will contribute $\omega(e)^2$, unless it is
  adjacent to a vertex of degree $0$, or unless it gets chosen later --- in
  each of these cases it contributes only $\omega(e)$.

  It follows that all inverse images of a marked floor diagram are tropical
  curves of the same multiplicity $\mult(C)$, and if there are $r$ inverse
  images we have $\mult(C)\cdot r=\mu(\D)\mu(\C)$.
\end{proof}

\begin {example} \label {ex-stirling}
  The following picture illustrates how we can complete a vertex of a marked
  Psi-floor diagram to floors of a tropical curve. The local picture of
  $\tilde{\D}$ on the left shows $2$ chosen incoming edges and $3$ non-chosen
  outgoing edges adjacent to a floor of degree $2$. The local multiplicity of
  this edge choice equals $\frac{2^0}{2!}\cdot \frac{2^3}{2!}=\frac{1}{2!}\cdot
  4$. We would like to complete this picture to the floor of a tropical curve.
  The lower part is unique. The factor of $\frac{1}{2!}$ for the lower part
  takes care of the two down ends which are adjacent to the contracted end and
  thus lead to a contribution of $\frac{1}{2!}$ in the factor $\nu_C$. For the
  upper part there are several possibilities. The middle column shows the
  $S(3,2)=3$ possibilities for $g=2$, i.e.\ for the case where all components
  obtained after removing the contracted marked edge also contain horizontal
  edges. The right column shows the $S(3,1)=1$ possibility for $g=1$, i.e.\ for
  the case where one of the ends of direction $(1,1)$ is directly adjacent to
  the contracted end. 

  \begin {center} \input {pics/stirling} \end {center}
\end {example}

\begin {remark}
  It follows immediately that also $\tilde N^{\floor}_{d,\kk} = \tilde
  N^{\trop}_{d,\kk}$ by taking the order of the contracted ends resp.\ vertices
  into account, both for the tropical curves and the floor diagrams.
\end {remark}

\subsection {Relative Psi-floor diagrams}

We now define relative analogues of Psi-floor diagrams and their markings. Fix
two sequences $\alpha$ and $\beta$. Our notation, which is more convenient for
our purposes, differs from \cite{FM}, where relative floor diagrams and their
markings were defined relative to partitions $\lambda = (1^{\alpha_1}
2^{\alpha_2} \cdots)$ and $\rho = (1^{\beta_1} 2^{\beta_2} \cdots)$. 

Let $\D$ be a Psi-floor diagram of degree $ d =  I(\alpha + \beta) $. A pair
$( \{ \alpha(v) \}, \{ \beta(v) \} )$ of collections of sequences, where $v$
runs over the vertices of $\D$, is called \emph {compatible} with $\D$ and
$(\alpha, \beta)$, if it satisfies:
\begin{enumerate}
\item The sums over each collection satisfy $\sum_{v \in V} \alpha(v) = \alpha$
  and $\sum_{v \in V} \beta(v) = \beta$.
\item For all vertices $v$ of $\D$ it holds that $I(\alpha(v)+\beta(v))= d_v -
  \dive(v)$.
\item If $d_v = 0$ then we require in addition that $|\alpha(v)| = 0$ and $
  | \beta(v) | = a_v + 2 -\val(v)$.
\end{enumerate}
The sequences $ \alpha(v) $ and $ \beta(v) $ correspond to the left (fixed and
non-fixed) ends adjacent to each floor. For a vertex of degree $0$, all
adjacent edges are directly adjacent to the contracted end, and thus there
cannot be any fixed ends in this case.

In the non-relative case, i.e.\ when $ \alpha = () $ and $ \beta = (d)$, it
necessarily follows that $ \alpha(v) = () $ and $ \beta(v) = (1 - \dive(v)) $
for all vertices $v$ of $\D$.

The \emph{(relative) type} $ \kk(\D)=(\kk_0,\kk_1,\ldots)$ of a Psi-floor
diagram $\D$ is defined as follows: for all $i\geq 1$ let $\kk_i$ be the number
of vertices $v$ of $\D$ with $a_v=i$. Set $\kk_0$ to be the number of vertices
with $a_v=0$ plus $2d+|\beta|-1-I\kk-\#V$. The latter number equals the number
of white vertices that we will add. This makes the equalities $ |\kk| =
2d+|\beta|-1-I\kk $, resp.\ $I(\alpha+\beta+\kk)=3d-1+|\beta|-|\kk|$ hold,
where the latter is equivalent to the former since $ d=I(\alpha+\beta) $.

The \emph{relative multiplicity} of a Psi-floor diagram $\D$ together with a
collection of sequences $\{ \beta(v) \}$ is
\begin {equation} \label {eq:relativediagrammult}
  \mu^{\rel} (\D)
   = \mu^{\rel} (\D,\{ \beta(v)\})
  := I^\beta
     \cdot \prod_{\text{edges }e} \omega(e)^2
     \cdot \prod_{\substack {
       v \, \stackrel e\to \, w \\[0.5ex]
       \text{s.t.\ } d_v = 0 \\[0.5ex]
       \text{or } d_w = 0
     }} \frac{1}{\omega(e)} \prod_{v:\, d_v = 0} \frac{1}{\beta(v)!}.
\end{equation}

For a collection of sequences $\{ \beta(v) \}$ and a vertex $v$ of $\D$ we
define the sets $I^{\text{rel}}(v)$ and $O^{\text{rel}}(v)$ by
\begin{displaymath}
  \begin{split}
    I^{\rel}(v) & := \{w \to v : d_w > 0 \}, \\
    O^{\rel}(v) & := \{v \to w: d_w > 0 \} \cup
      \coprod \, \{ v \stackrel{i}{\to} \circ \},
  \end{split}
\end{displaymath}
where the latter is a disjoint union of the outgoing edges of $\D$ at $v$
augmented by $\beta_i^v$ \emph{indistinguishable} edges of weight $i$ for all $
i \ge 1 $, directed away from $v$ and ending in distinct vertices $\circ$.
These indistinguishable edges correspond to the non-fixed ends of the tropical
curve adjacent to a floor, which a priori could be adjacent to the contracted
end, and therefore can be chosen.

\begin{example} \label{ex:relativeinandout}
  Below we have indicated the sets $I^\text{rel}(v)$ and $O^\text{rel}(v)$ in
  the case of the Psi-floor diagram of Example \ref{ex:floordiagram} with
  $\alpha = (1)$, $\beta=(2,1)$, and all $\alpha(v)$ and $\beta(v)$ being the
  zero sequence unless indicated otherwise. The relative multiplicity $
  \mu^{\rel} (\D,\{\beta(v)\}) $ is $ 4 \cdot 2 = 8 $.

  \begin{displaymath}
    \begin{picture}(50,75)(30,-50)\setlength{\unitlength}{5pt}\thicklines
    \oooo\Eeee\eEee\eeEe
    \put(15,1.5){\makebox(0,0){$2$}} 
    \put(5,0){\vector(1,0){1}} 
    \put(15,0){\vector(1,0){1}} 
    \put(25,0){\vector(1,0){1}} 
    \put(-9.4,-3){\makebox(0,0){$(d_v, a_v)= $}}
    \put(0,-3){\makebox(0,0){$1 \, 0$}}
    \put(10,-3){\makebox(0,0){$2 \, 3$}}
    \put(20,-3){\makebox(0,0){$1 \, 2$}}
    \put(30,-3){\makebox(0,0){$1 \, 0$}}
    \put(-8,-6){\makebox(0,0){$\alpha(v)=$}}
    \put(30,-6){\makebox(0,0){$(1)$}}
    \put(-8,-9){\makebox(0,0){$\beta(v)=$}}
    \put(10,-9){\makebox(0,0){$(1)$}}
    \put(20,-9){\makebox(0,0){$(0,1)$}}
    \put(30,-9){\makebox(0,0){$(1)$}}
    \put(10,+0.0){\line(1,1){4.2}}
    \put(10.7,+0.7){\vector(1,+1){3}}
    \put(15,+5){\circle{2}}
    \put(20,+0.0){\line(1,1){4.2}}
    \put(20.7,+0.7){\vector(1,+1){3}}
    \put(25,+5){\circle{2}}
    \put(21.7,4){\makebox(0,0){$2$}}
    \put(30,+0.0){\line(1,1){4.2}}
    \put(30.7,+0.7){\vector(1,+1){3}}
    \put(35,+5){\circle{2}}
    \end{picture}
  \end{displaymath}
\end{example}

As before, an \emph{edge choice} $\C(\D)$ is given by a subset $C(v) \subset
I^{\text{rel}}(v) \cup O^{\text{rel}}(v)$ for each floor $v$ of $\D$ such that
$|C(v)| = a_v +2 - 2 d_v $ for all $v$, and $C(v) \cap C(w) = \emptyset$ for
distinct floors $v$ and $w$. If $d_v = 0$, we set $C(v) = \emptyset$. The
\emph{local multiplicity at} $v$ of such a choice is
\begin{equation} \label {eq:localmult}
  \mu^{\text{rel}}_{v, C(v)} := \begin {cases}
    \frac {d_v^{i(v)}}{d_v!} \cdot \frac{d_v^{o(v)}}{d_v!}
      & \mbox {if $ d_v > 0 $}, \\[0.5ex]
    1 & \mbox {if $ d_v = 0 $},
  \end{cases}
\end{equation}
where, similarly to the absolute case, $i(v) = |I^{\text{rel}}(v) \backslash
C(v)|$ is the number of non-chosen incoming edges and $o(v) =
|O^{\text{rel}}(v) \backslash C(v)| + |\alpha(v)|$ is the number of non-chosen
edges in $O^{\text{rel}}(v)$ together with some additional edges (corresponding
to tangency conditions at fixed points, resp.\ to fixed left ends).

The \emph{relative multiplicity} of the edge choice $\C$ of the Psi-floor
diagram $\D$ together with a compatible pair of collections of sequences $(\{
\alpha(v) \}, \{ \beta(v)\})$ is
\begin{equation} \label{eq:relativechoicemult}
  \mu^{\rel}(\C)
    := \mu^{\rel}(\C, \{ \alpha(v) \}, \{ \beta(v) \})
    := \prod_{v \in V} \mu_{v, C(v)}^{\rel} \,
       \prod_{v \in V} \prod_{e \in C(v)} \frac{1}{\omega(e)} \,
       \prod_{v \in V} \frac{1}{c(v) !},
\end{equation}
where $c(v)$ is the sequence given by $c(v)_i := |C(v) \cap \{v \stackrel{i}{
\to} \circ \}|$ for $i \ge 1$.

{\bf Example \ref{ex:relativeinandout} (continued).}
  An example of an edge choice for the above Psi-floor diagram together with
  collections $\{ \alpha(v) \}$ and $\{ \beta(v) \}$ is given below. As before,
  we indicate chosen edges by thickening edges at the vertices where they are
  chosen. Notice that $|C(v)| = a_v - 2(d_v-1)$ at every vertex $v$ since there
  are no vertices of degree $0$. The relative multiplicity of the edge choice
  is $ \mu^{\rel} (\C) = \tfrac{1}{2}$.

  \begin{displaymath}
    \begin{picture}(50,30)(65,-15)\setlength{\unitlength}{5pt}\thicklines
    \oooo\Eeee\eEee\eeEe
    \put(15,1.5){\makebox(0,0){$2$}} 
    \put(5,0){\vector(1,0){1}} 
    \put(15,0){\vector(1,0){1}} 
    \put(25,0){\vector(1,0){1}} 
    \put(0,-3){\makebox(0,0){$1 \, 0$}}
    \put(10,-3){\makebox(0,0){$2 \, 3$}}
    \put(20,-3){\makebox(0,0){$1 \, 2$}}
    \put(30,-3){\makebox(0,0){$1 \, 0$}}
    \put(10,+0.0){\line(1,1){4.2}}
    \put(10.7,+0.7){\vector(1,+1){3}}
    \put(15,+5){\circle{2}}
    \put(20,+0.0){\line(1,1){4.2}}
    \put(20.7,+0.7){\vector(1,+1){3}}
    \put(25,+5){\circle{2}}
    \put(22,4){\makebox(0,0){$2$}}
    \put(30,+0.0){\line(1,1){4.2}}
    \put(30.7,+0.7){\vector(1,+1){3}}
    \put(35,+5){\circle{2}}
    \linethickness{1mm}
    \put(6,0){\line(1,0){4}}
    \put(16,0){\line(1,0){4}}
    \put(20,0){\qbezier(0,0)(1,1)(2.5,2.5)}
    \end{picture}
  \end{displaymath}

\begin {definition} \label{def:relativemarking}
  An \emph{$(\alpha, \beta)$-marking} of a Psi-floor diagram $\D$ with a
  compatible choice of a pair of collections $(\{\alpha(v)\}, \{\beta(v)\})$
  and an edge choice $\C(\D)$ is defined by the following three-step process
  which we illustrate in the case of Example \ref{ex:relativeinandout}.

  {\bf Step 1:} For each vertex $v$ of $\D$ and every $i \ge 1$ create
  $\beta(v)_i - |C(v) \cap \{ v \stackrel{i}{\to} \circ \}|$ new vertices
  (which we call \emph{$\beta$-vertices} and illustrate as
  \begin{picture}(10,8)(0,0) \thicklines
    \put(5,4){\circle{8}}
  \end{picture}),
  and connect them to $v$ with new edges of weight $i$ directed away from $v$.
  Similarly, create $\alpha(v)_i $ new vertices (which we call
  \emph{$\alpha$-vertices} and illustrate as
  \begin{picture}(10,8)(0,0) \thicklines
    \put(5,4){\circle{8}}
    \put(5,4){\circle*{4}}
  \end{picture})
  and connect them to $v$ with new edges of weight $i$ directed away from $v$.

  \begin{displaymath}
    \begin{picture}(50,40)(67.5,-20)\setlength{\unitlength}{5pt}\thicklines
    \oooo\Eeee\eEee\eeEe
    \put(15,1.5){\makebox(0,0){$2$}} 
    \put(5,0){\vector(1,0){1}} 
    \put(15,0){\vector(1,0){1}} 
    \put(25,0){\vector(1,0){1}} 
    \put(0,-3){\makebox(0,0){$1 \, 0$}}
    \put(10,-3){\makebox(0,0){$2 \, 3$}}
    \put(20,-3){\makebox(0,0){$1 \, 2$}}
    \put(30,-3){\makebox(0,0){$1 \, 0$}}
    \put(10,+0.0){\line(1,1){4.2}}
    \put(10.7,+0.7){\vector(1,+1){3}}
    \put(15,+5){\circle{2}}
    \put(30,+0.0){\line(1,1){4.2}}
    \put(30.7,+0.7){\vector(1,+1){3}}
    \put(35,+5){\circle{2}}
    \put(30,+0.0){\line(2,1){9}}
    \put(30,+0){\vector(2,+1){7}}
    \put(40,+5){\circle{2}}
    \put(40,+5){\circle*{1}}
    \linethickness{1mm}
    \put(6,0){\line(1,0){4}}
    \put(16,0){\line(1,0){4}}
    \end{picture}
  \end{displaymath}

  {\bf Step 2:} Subdivide each non-chosen edge of the original Psi-floor
  diagram $\D$ between floors into two edges by introducing a new vertex for
  each edge. The new edges inherit their weights and orientations. Call the
  resulting graph $\tilde{\D}$.

  \begin{displaymath}
    \begin{picture}(50,40)(67.5,-20)\setlength{\unitlength}{5pt}\thicklines
    \oooo\Eeee\eEee
    \put(15,1.5){\makebox(0,0){$2$}} 
    \put(5,0){\vector(1,0){1}} 
    \put(15,0){\vector(1,0){1}} 
    \put(20,+0.0){\line(1,0){4}}
    \put(23,0){\vector(1,0){1}} 
    \put(25,0){\circle{2}}
    \put(26,+0.0){\line(1,0){4}}
    \put(28,0){\vector(1,0){1}} 
    \put(0,-3){\makebox(0,0){$1 \, 0$}}
    \put(10,-3){\makebox(0,0){$2 \, 3$}}
    \put(20,-3){\makebox(0,0){$1 \, 2$}}
    \put(30,-3){\makebox(0,0){$1 \, 0$}}
    \put(10,+0.0){\line(1,1){4.2}}
    \put(10.7,+0.7){\vector(1,+1){3}}
    \put(15,+5){\circle{2}}
    \put(30,+0.0){\line(1,1){4.2}}
    \put(30.7,+0.7){\vector(1,+1){3}}
    \put(35,+5){\circle{2}}
    \put(30,+0.0){\line(2,1){9}}
    \put(30,+0){\vector(2,+1){7}}
    \put(40,+5){\circle{2}}
    \put(40,+5){\circle*{1}}
    \linethickness{1mm}
    \put(6,0){\line(1,0){4}}
    \put(16,0){\line(1,0){4}}
    \end{picture}
  \end{displaymath}

  {\bf Step 3:} Order the vertices of $\tilde{\D}$ linearly, extending the
  order of the vertices of the original Psi-floor diagram $\D$, such that (as
  in $\D$) each edge is directed from a smaller vertex to a larger vertex.
  Furthermore, we require that the $\alpha$-vertices are largest among all
  vertices, and for every pair of $\alpha$-vertices $ v>w $ the weight of the
  $v$-adjacent edge is larger than or equal to the weight of the $w$-adjacent
  edge.

  \begin{displaymath}
    \begin{picture}(50,40)(60,-20)\setlength{\unitlength}{5pt}\thicklines
    \put(0,+0){\circle*{2}}
    \put(5,+0){\circle*{2}}
    \put(10,+0){\circle*{2}}
    \put(15,+0){\circle{2}}
    \put(20,+0){\circle{2}}
    \put(25,+0){\circle*{2}}
    \put(30,+0){\circle{2}}
    \put(35,+0){\circle{2}}
    \put(35,+0){\circle*{1}}
    \put(0,-3){\makebox(0,0){$1 \, 0$}}
    \put(5,-3){\makebox(0,0){$2 \, 3$}}
    \put(10,-3){\makebox(0,0){$1 \, 2$}}
    \put(25,-3){\makebox(0,0){$1 \, 0$}}
    \put(0,+0.0){\line(1,0){5}}
    \put(2,0){\vector(1,0){1}}
    \put(5,+0.0){\line(1,0){4}}
    \put(7,0){\vector(1,0){1}}
    \put(7.5,1.2){\makebox(0,0){$2$}}
    \put(21,+0.0){\line(1,0){4}}
    \put(23,0){\vector(1,0){1}}
    \put(25,+0.0){\line(1,0){4}}
    \put(28,0){\vector(1,0){1}}
    \qbezier(5.8,0.6)(7,3)(10,3)\qbezier(10,3)(13,3)(14.2,0.6)
    \put(10,3){\vector(1,0){1}}
    \qbezier(10.8,-0.6)(12,-3)(15,-3)\qbezier(15,-3)(18,-3)(19.2,-0.6)
    \put(15,-3){\vector(1,0){1}}
    \qbezier(25.8,0.6)(27,3)(30,3)\qbezier(30,3)(33,3)(34.2,0.6)
    \put(30,3){\vector(1,0){1}}
    \linethickness{1mm}
    \put(3,0){\line(1,0){2}}
    \put(8,0){\line(1,0){2}}
    \end{picture}
  \end{displaymath}

  The (in this example unique) tropical curve mapping to the floor diagram
  above can be found in Example \ref{ex-reltropcurve}. As in the non-relative
  case, we call the extended graph $\tilde{\D}$ together with the linear order
  on its vertices an \emph{$(\alpha, \beta)$-marked Psi-floor diagram}, or an
  \emph{$(\alpha, \beta)$-marking} of the Psi-floor diagram $\D$.
\end{definition}

In step 1 we added $|\beta|$ white vertices (of which we later remove the
chosen ones), and in step 2 we subdivide the non-chosen ones of the $\#V-1$
bounded edges. That is, altogether we added $|\beta|+\#V-1-\sum_{v\in
V}(a_v-2(d_v-1))=2d-1+|\beta|-I\kk-\#V$ white vertices.

As before, we need to count $(\alpha,\beta)$-marked Psi-floor diagrams up to
equivalence. Two $(\alpha,\beta)$-marked Psi-floor diagrams $\tilde{\D}_1$,
$\tilde{\D}_2$ are \emph{equivalent} if $\tilde{\D}_1$ can be obtained from
$\tilde{\D}_2$ by permuting edges without changing their weights, i.e.\ if
there exists an automorphism of weighted graphs which preserves the vertices of
$\D$ and maps $\tilde{\D}_1$ to $\tilde{\D}_2$. The \emph{number of markings}
$\nu^\text{rel}(\D, \C) = \nu^\text{rel}(\D, \{ \alpha(v) \},\{ \beta(v) \},
\C)$ is the number of $(\alpha,\beta)$-marked Psi-floor diagrams $\tilde{\D}$
up to equivalence. In our running example we have $ \nu^{\rel}(\D,\C) = 5 $:
the white vertex attached to the floor labeled $ (2,3) $ can be placed in the
linear order at any position to the right of this floor and to the left of the
$ \alpha $-vertex.

By specializing to the case $a_v = 0$ for all vertices $v$ of $\D$ we recover
the definition of $(\lambda, \rho)$-markings of floor diagrams of Fomin and
Mikhalkin \cite{FM}, for partitions $\lambda = (1^{\alpha_1} 2^{\alpha_2}
\cdots)$ and $\rho = (1^{\beta_1} 2^{\beta_2} \cdots)$. As in the non-relative
case, all floors necessarily have degree $d_v = 1$ and no edges get chosen.

\begin {definition}[$ N^{\floor}_{d,\kk}(\alpha,\beta) $ and $
    \tilde N^{\floor}_{d,\kk}(\alpha,\beta) $]
  Let $d \ge 1$ and $\alpha, \beta$ be two sequences with $ I(\alpha+\beta) =
  d $. Furthermore, let ${\bf k}$ be a sequence of non-negative integers with
  $I(\alpha+\beta+\kk)= 3d -1 +|\beta|-|\kk| $. Set
  \begin {displaymath}
    N^{\floor}_{d, {\bf k}}(\alpha, \beta)
      := \sum_{\D, \{ \alpha(v) \}, \{ \beta(v) \}} \mu^{\text{rel}}(\D) \,
         \sum_\C \mu^\text{rel}(\C) \, \nu^\text{rel}(\D, \C),
  \end {displaymath}
  where the first sum is over all degree $d$ Psi-floor diagrams of type ${\bf
  k}$ and over all compatible pairs of collections $ (\{\alpha(v)\},
  \{\beta(v)\})$, and the second sum is over all edge choices $\C$ of $\D$.
  Correspondingly (see Definition \ref {def-ndk-ab}), we set $ \tilde
  N^{\floor}_{d,\kk}(\alpha,\beta) := \beta! \cdot \frac {\kk!}{|\kk|\,!} \,
  N^{\floor}_{d,\kk}(\alpha,\beta) $.
\end {definition}

\begin {remark} \label {rem-relnfloortilde}
  As in Remark \ref{rem-nfloortilde}, we can also define the numbers
  $\tilde{N}^{\floor}_{d, {\bf k}}(\alpha,\beta)$ directly using Psi-floor
  diagrams. Then we require that Psi-powers of the vertices of the marked
  Psi-floor diagram are in a fixed order, and we mark the white end vertices.
\end {remark}

\begin {theorem}[The equality $ N^{\floor}_{d,\kk}(\alpha, \beta) =
    N^{\trop}_{d,\kk}(\alpha, \beta)$] \label {thm:relativepsicount}
  Let $d \ge 1$ and $\alpha, \beta$ be two sequences with $I(\alpha+\beta) =
  d$. Let ${\bf k}$ be a sequence of non-negative integers satisfying
  $I(\alpha+\beta+\kk)=3d-1+|\beta|-|\kk|$. Then $ N^{\floor}_{d, {\bf k}}(
  \alpha, \beta)=  N^{\trop}_{d, {\bf k}}(\alpha, \beta)$.
\end {theorem}

The proof is analogous to the proof of Theorem \ref{thm:psicount}.

\begin {remark}
  Again, it follows immediately that the same equality holds for the numbers
  $\tilde{N}^{\floor}_{d, {\bf k}}(\alpha,\beta)=\tilde{N}^{\trop}_{d, {\bf
  k}}(\alpha,\beta)$ as well.
\end {remark}

\subsection {The Caporaso-Harris formula for floor diagrams}
  \label {sec-CHfloor}

Now we use Psi-floor diagrams to obtain the Caporaso-Harris type recursion of
Corollary \ref{cor-CHclass} for the numbers $N_{d, \kk}^{\floor}(\alpha,
\beta)$. As this recursion formula determines all the numbers it follows that $
N_{d,\kk}(\alpha, \beta) = N_{d,\kk}^{\floor}(\alpha, \beta)$. As we know by
Theorem \ref{thm:relativepsicount} that also $ N_{d,\kk}^{\floor}(\alpha,
\beta) = N_{d,\kk}^{\trop}(\alpha, \beta)$ holds, we thus have that
  \[ N_{d,\kk}(\alpha, \beta)
     = N_{d,\kk}^{\floor}(\alpha, \beta)
     = N_{d,\kk}^{\trop}(\alpha, \beta) \]
for all $ d,\kk,\alpha,\beta $, as claimed in Remark \ref {rem-equiv-trop-alg}.
We use Notation \ref {not-sequences} and the notation in equation (\ref
{eq-binoms}) below.

\begin {theorem}[Caporaso-Harris formula for Psi-floor diagrams]
    \label{thm-CHfloor}
  The numbers $N^{\floor}_{d, {\bf k}}(\alpha, \beta)$ satisfy the
  Caporaso-Harris recursion in Corollary \ref{cor-CHclass}.
\end {theorem}

\begin {proof}
  The basic strategy is to examine the possibilities for the largest vertex
  $v'$ of an $(\alpha, \beta)$-marking $\tilde{\D}$ of a Psi-floor diagram $\D$
  of degree $d$ and type ${\bf k}$ which is not an $\alpha$-vertex (see step 1
  in Definition \ref{def:relativemarking} to recall the definition of
  $\alpha$-vertices and $\beta$-vertices). The idea is to ``cut off'' the
  vertex $v'$ and to interpret the contributions of the connected components of
  the remaining part again in terms of smaller floor diagrams. 

  The complement of $v'$ and the $v'$-adjacent edges in $\tilde{\D}$ consists
  of markings $\tilde{\D^1}, \dots, \tilde{\D^t}$ of Psi-floor diagrams $\D^1,
  \dots, \D^t$ and some isolated $\alpha$-vertices. For $1 \le i \le t$ define
  \begin{enumerate}
  \item $d^i$ and ${\bf k}^i$ to be the degree and the type of $\D^i$,
    respectively,
  \item $ \alpha^i = \sum \alpha(v) $ to be the sequence of multiplicities of
    edge weights between $\D^i$ and the $\alpha$-vertices of $\tilde{\D}$,
    where the sum is over all vertices $v$ in the Psi-floor diagram $\D^i$, 
  \item $ \beta^i = \sum \beta(v) $, the respective count for the
    $\beta$-vertices of $\tilde{\D}$,
  \item $m^i$ to be the weight of the edge between $v'$ and $\D^i$.
  \end{enumerate}
  Of course, $m^i= d^i-I(\alpha^i+\beta^i)$.

  We will see later that all contributions from the components $ \D^i $ are of
  the form $ N^{\floor}_{d^i,\bf{k^i}}(\alpha^i + e_{m^i}, \beta^i)$ resp.\ $
  N^{\floor}_{d^i, \bf{k^i}}(\alpha^i, \beta^i  + e_{m^i}) $. In these cases we
  necessarily have 
  \begin {equation} \label{eq:dimensionscount1}
    I\alpha^i+m^i+I\beta^i+I\kk^i= 3d^i-1+|\beta^i|- |\kk^i|, \mbox{ resp.}
  \end {equation}
  \begin {equation} \label{eq:dimensionscount2}
    I\alpha^i+I\beta^i+m^i+I\kk^i= 3d^i-1+|\beta^i|+1- |\kk^i|.
  \end {equation}

  Now consider the possibilities for the largest vertex $v'$. We will
  distinguish three cases.
  
  \textbf {Case 1:} The vertex $v'$ is not a vertex of the original diagram
  $\D$. Hence $ \tilde\D $ looks locally around $ v' $ as in the following
  picture.

  \begin {center} \input {pics/CH1} \end {center}

  Then $t = 1$, $\alpha^1 = \alpha$ and $\beta^1 = \beta - e_{m^1}$. The
  $(\alpha^1,\beta^1)$-markings of $\D$ with $v'$ maximal among all
  non-$\alpha$-vertices are in canonical bijection with $(\alpha^1 + e_{m^1},
  \beta^1)$-markings of $\D$ (by making $v'$ an $\alpha$-vertex and, for
  example, inserting it to the right of the other $\alpha$-vertices adjacent to
  weight $m^1$ edges). This bijection is weight-preserving up to a factor
  $m^1$, as edges of weight $m^1$ adjacent to $\beta$-vertices contribute a
  factor of $m^1$ whereas edges adjacent to $\alpha$-vertices do not (see
  equation (\ref{eq:relativediagrammult})). Thus, if $v'$ is not a vertex of
  the original diagram we get a contribution of
    \[ \sum_{m^1: \, \beta_{m^1} > 0} m^1 \cdot N^{\floor}_{d^1,
       {\bf k}^1}(\alpha^1 + e_{m^1}, \beta^1). \]
  This contribution equals the summands with $d'=0$ and $a=0$ in the sum of
  Corollary \ref{cor-CHclass}: for $d' = 0$ the non-vanishing of
  $d'^{|\alpha'|+t-t'}$ implies that $|\alpha'|=0$ and $t=t'$, and
  equation (\ref{eq:dimensionscount1}) (which can be rearranged to imply a
  valence and divergence condition on $v'$ as we will show below) implies
  furthermore that $t'=1$. This finishes case 1.

  Now assume that $v'$ is a vertex of the original diagram $\D$, and denote by
  $d'$ and $a$ the degree and Psi-power of $v'$, respectively. We need to
  count the number of ways in which markings of the Psi-floor diagrams $\D^1,
  \dots,\D^t$ can be combined to a marking of the Psi-floor diagram $\D$. We
  need to distinguish whether $v'$ is a floor of $\D$ (i.e.\ $d' > 0$) or not.

  \textbf {Case 2:} $v'$ is a vertex of $\D$, and $ d'=0 $. Then we obtain the
  following local picture for $ \tilde\D $.

  \begin {center} \input {pics/CH2} \end {center}

  In this case none of the edges between $v'$ and the Psi-floor diagrams $\D^i$
  can be chosen. Notice that $(\alpha^i,\beta^i+e_{m^i})$-markings of $\D^i$
  with $v'$ largest among all $\beta$-vertices (if we consider $v'$ as a
  $\beta$-vertex of $\D^i$) are in canonical bijection with $(\alpha^i +
  e_{m^i}, \beta^i)$-markings of $ \D^i $. This bijection is weight-preserving
  up to a factor of $m^i$ (see equation (\ref{eq:relativediagrammult})).

  To count the number of ways in which we can combine the markings of the
  pieces fix an $(\alpha^i + e_{m^i}, \beta^i)$-marking of $\D^i$, one for each
  $1 \le i \le t$. Produce an $(\alpha,\beta)$-marking of $\D$ as follows:
  First, glue the markings by identifying all largest $\alpha$-vertices in each
  of the marking of $\D^i$ adjacent to an edge of weight $m^i$ with each other
  (thereby obtaining the vertex $v'$). Then order the $\alpha$-vertices of the
  markings by extending the partial order on the set of $\alpha$-vertices given
  by the markings of the components to a linear order on all vertices. There
  are $\binom{\alpha}{\alpha^1, \dots, \alpha^t}$ ways to do this.

  In a second step, we extend the partial order on the vertices that are less
  than $v'$ to a linear order on all vertices less than $v'$. As $v'$ is
  maximal among the non-$\alpha$-vertices of the marking $\tilde{\D}$ it has
  $|{\bf k}| - 1$ vertices which are less than $v'$. Using the earlier
  bijection between $(\alpha^i, \beta^i)$-markings of $\D^i$ with $v'$ largest
  among all $\beta$-vertices (if we consider $v'$ as a $\beta$-vertex of
  $\D^i$) and $(\alpha^i + e_{m^i}, \beta^i)$-markings of $ \D^i $ we see
  that there are  $|{\bf k}^i|$ vertices in component $i$ which are less than
  $v'$. Hence there are $\binom {|{\bf k}| - 1}{|{\bf k}^1|, \dots, |{\bf
  k}^t|}$ linear extensions of the partial order that is induced by the
  linear orders of the components.

  By equation (\ref{eq:relativediagrammult}) the product of the contributions
  from the $t$ components differs from the contribution of the marking $
  \tilde{\D}$ by $\frac{1}{\beta(v')!}$, but $\beta(v') = \beta - \sum
  \beta^i=\beta'$. Moreover, we overcount by $t!$ as we labeled the unlabeled
  components $1, \dots, t$. Altogether, we get a contribution of
    \[ \sum \frac{1}{t!} \, \frac{1}{\beta' !} \,
            \binom {|{\bf k}| - 1}{|{\bf k}^1|, \dots, |{\bf k}^t|} \cdot
            \binom{\alpha}{\alpha^1, \dots, \alpha^t} \,
            \prod_{i=1}^t m^i \,
            \prod_{i=1}^t N^{\floor}_{d^i,\kk^i}(\alpha^i+e_{m^i},\beta^i), \]
  which equals the summands with $d'=0$ but $a>0$ in the recursion of Corollary
  \ref{cor-CHclass}. As before, equations (\ref{eq:dimensionscount1}) and
  (\ref{eq:dimensionscount2}) imply that $v'$ has the correct divergence and
  valence (see below).

  \textbf {Case 3:} $v'$ is a vertex of $\D$, and $ d'>0 $. In this case we
  obtain the following local picture for $ \tilde\D $.

  \begin {center} \input {pics/CH3} \end {center}

  As then $v'$ is largest among all non-$\alpha$-vertices we have $C(v')
  \supset O(v')$. Define the Psi-floor diagrams $\D^1, \dots, \D^t$ and their
  markings $\tilde{\D^1}, \dots, \tilde{\D^t}$ as before, as well as $d^i$,
  ${\bf k}^i$, $\alpha^i$, $\beta^i$ and $m^i$, for $1 \le i \le t$. Without
  loss of generality we can assume that there is a number $ t' \in \{0,\dots,t
  \} $ such that the edges between $v'$ and $\tilde{\D^i}$ are chosen at $v'$
  for all $ i \le t' $, whereas for $ i > t' $ they are not.

  Now consider a component $ \D^i $. We treat the cases $i \le t'$ and $ i>t' $
  separately. If $i \le t'$ then the $(\alpha^i, \beta^i + e_{m^i})$-markings
  of $\D^i$ with $v'$ largest among all $\beta$-vertices (if we consider $v'$
  as a $\beta$-vertex of $\D^i$) are in canonical bijection with $(\alpha^i +
  e_{m^i}, \beta^i)$-markings of $\D^i$ by the same reasoning as for $d' = 0$.
  As we have seen, this bijection is weight-preserving up to a factor of $m^i$.

  If $i > t'$ then a linear order (up to equivalence) on the vertices of
  $\tilde{\D^i}$ that can be extended to a marking of $\D$ canonically
  determines an $(\alpha^i, \beta^i + e_{m^i})$-marking of $\D^i$ together
  with a distinguished $\beta$-vertex adjacent to an edge of weight $m^i$
  (namely the image of the edge of $\tilde{\D^i}$ which is closest to $v'$ in
  $\tilde{\D}$). Conversely, given an $(\alpha^i, \beta^i + e_{m^i})$-marking
  of $\D^i$ together with a distinguished $\beta$-vertex adjacent to an edge of
  weight $m^i$, this canonically determines a linear order (up to equivalence)
  on the vertices of $\tilde{\D^i}$ that can be extended to a marking of $\D$.
  This $(\beta_{m^i}^i + 1)$-to-$1$ map is weight-preserving up to a factor of
  $m^i$.

  Again, to produce an $(\alpha, \beta)$-marking of $\D$ we need to extend the
  partial order on the set of $\alpha$-vertices given by the markings of the
  components to a linear order on all $\alpha$-vertices. There is no difference
  to the $d' = 0$ case, hence there are $\binom{\alpha}{\alpha^1, \dots,
  \alpha^t}$ different extensions. As before, there are $\binom{|\kk| -
  1}{\kk^1, \dots, \kk^t}$ ways to extend the partial order on the vertices
  that are smaller than $v'$ to a linear order.

  Also as before, by equation (\ref{eq:relativediagrammult}) the weight of a
  marking of $\D$ differs from the product of the individual weights of the $t$
  components by contributions from the vertex $v'$. The local multiplicity at
  $v'$ from equation (\ref{eq:localmult}) is $\frac{(d')^{t - t'}}{d' !}
  \frac{(d')^{|\alpha'|}}{d'!} $ as the number of non-chosen incoming vertices
  is $i(v') = |I(v') \backslash C(v')| = t - t'$ and 
  \begin{displaymath}
    o(v') = |O(v') \backslash C(v')| + |\alpha(v')|
          = 0 + |\alpha'| = |\alpha'|
  \end{displaymath}
  since $\alpha(v') = \alpha - \sum_i \alpha^i $. The second contribution
  from the vertex $v'$ is $\frac{1}{\beta' !}$ (see equation (\ref
  {eq:relativechoicemult})), as $\beta' = \beta(v')$, $C(v') \supset O(v')$,
  and hence $c(v') = \beta(v')$, and these are the only contributions in which
  the markings of $\D$ and the contributions from its components differ.
  Moreover, we overcount by $t'! \cdot (t-t')!$ as we labeled the unlabeled
  components $1, \dots, t'$ and $t'+1, \dots, t$.

  \textbf {Divergence and valence conditions for all cases:} In all three
  cases, equations (\ref {eq:dimensionscount1}) and (\ref
  {eq:dimensionscount2}) imply that $v'$ has the correct divergence and
  valence: summing up equations (\ref {eq:dimensionscount1}) for $1 \le i \le
  t'$ and (\ref {eq:dimensionscount2}) for $t' + 1 \le i \le t$ yields
    \[ I\alpha-I\alpha'+I\beta-I\beta'+I\kk -a +\sum m^i
         = 3d-3d'-t'+|\beta|-|\beta'|-|\kk|+1. \]
  Since $I(\alpha+\beta+\kk)=3d-1+|\beta|-|\kk|$ we can conclude
  \begin{align} \label{eq-divergence}
    -I\alpha'-I\beta' -a +\sum m^i = -3d'-t'-|\beta'|+2.
  \end{align}
  Now replace $m^i$ by $d^i-I(\alpha^i+\beta^i)$ and use that $d = I(\alpha +
  \beta) $ to obtain the valence condition at $v'$:
    \[ -a = -2d'-t'-|\beta'|+2, \mbox{ resp. } 2d'+t'+|\beta'|=a+2. \]
  Together with equation (\ref{eq-divergence}) the valence condition implies
  the divergence condition at $v'$:
    \[ -I\alpha'-I\beta'  +\sum m^i= -d', \mbox{ resp. } d'+\sum m^i
         = I(\alpha'+\beta'). \]
  Hence the contributions in the case when $v'$ is a floor equal the summands
  with $d'>0$ in the recursion of Corollary \ref{cor-CHclass}. This completes
  the proof.
\end {proof}

Of course, one can also prove the recursion in Theorem \ref {thm-CHclass}
directly using Psi-floor diagrams. We then have to use the numbers $
\tilde{N}^{\floor}_{d, {\bf k}}(\alpha,\beta) $ of Remark \ref
{rem-relnfloortilde}, where we fix an order for the Psi-powers and mark the
white end vertices. 

%% file: pics/floors.tex
\begin{picture}(0,0)%
\includegraphics{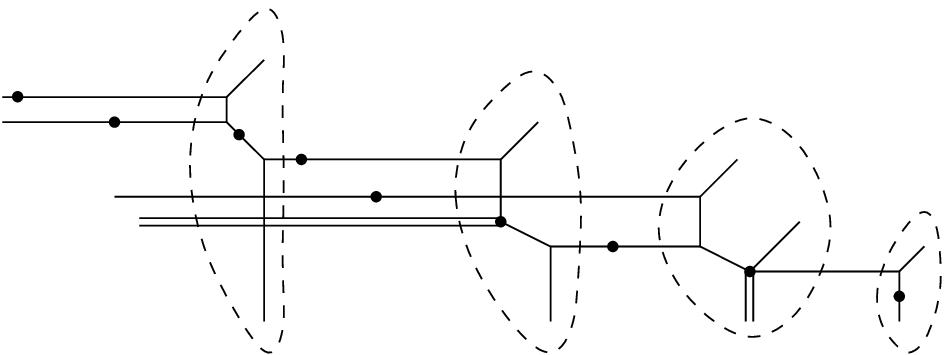}%
\end{picture}%
\setlength{\unitlength}{3947sp}%
\begingroup\makeatletter\ifx\SetFigFont\undefined%
\gdef\SetFigFont#1#2#3#4#5{%
  \reset@font\fontsize{#1}{#2pt}%
  \fontfamily{#3}\fontseries{#4}\fontshape{#5}%
  \selectfont}%
\fi\endgroup%
\begin{picture}(4528,1674)(3964,-6523)
\put(6969,-5966){\makebox(0,0)[lb]{\smash{{\SetFigFont{8}{9.6}{\familydefault}{\mddefault}{\updefault}{\color[rgb]{0,0,0}$2$}%
}}}}
\end{picture}%

%% file: pics/zerofloor.tex
\begin{picture}(0,0)%
\includegraphics{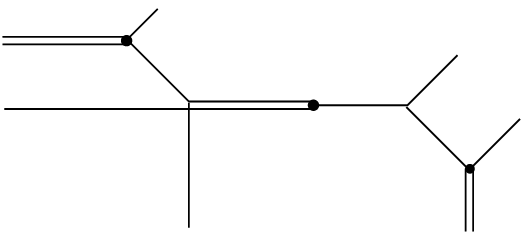}%
\end{picture}%
\setlength{\unitlength}{3947sp}%
\begingroup\makeatletter\ifx\SetFigFont\undefined%
\gdef\SetFigFont#1#2#3#4#5{%
  \reset@font\fontsize{#1}{#2pt}%
  \fontfamily{#3}\fontseries{#4}\fontshape{#5}%
  \selectfont}%
\fi\endgroup%
\begin{picture}(2507,1092)(4195,-5941)
\put(5835,-5288){\makebox(0,0)[lb]{\smash{{\SetFigFont{8}{9.6}{\familydefault}{\mddefault}{\updefault}{\color[rgb]{0,0,0}$2$}%
}}}}
\end{picture}%

%% file: pics/stirling.tex
\begin{picture}(0,0)%
\includegraphics{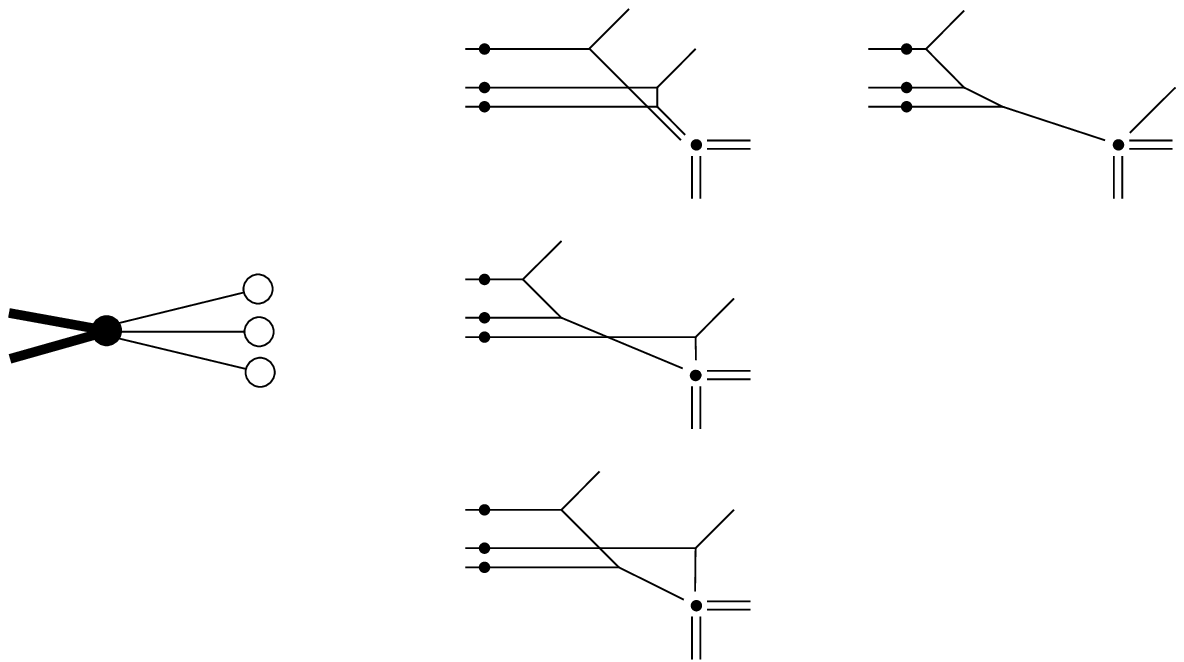}%
\end{picture}%
\setlength{\unitlength}{3947sp}%
\begingroup\makeatletter\ifx\SetFigFont\undefined%
\gdef\SetFigFont#1#2#3#4#5{%
  \reset@font\fontsize{#1}{#2pt}%
  \fontfamily{#3}\fontseries{#4}\fontshape{#5}%
  \selectfont}%
\fi\endgroup%
\begin{picture}(5654,3146)(369,-2979)
\put(4549,  3){\makebox(0,0)[lb]{\smash{{\SetFigFont{8}{9.6}{\familydefault}{\mddefault}{\updefault}{\color[rgb]{0,0,0}$2$}%
}}}}
\put(2602,-1094){\makebox(0,0)[lb]{\smash{{\SetFigFont{8}{9.6}{\familydefault}{\mddefault}{\updefault}{\color[rgb]{0,0,0}$2$}%
}}}}
\put(2604, 11){\makebox(0,0)[lb]{\smash{{\SetFigFont{8}{9.6}{\familydefault}{\mddefault}{\updefault}{\color[rgb]{0,0,0}$2$}%
}}}}
\put(2602,-2201){\makebox(0,0)[lb]{\smash{{\SetFigFont{8}{9.6}{\familydefault}{\mddefault}{\updefault}{\color[rgb]{0,0,0}$2$}%
}}}}
\put(884,-1640){\makebox(0,0)[b]{\smash{{\SetFigFont{10}{12.0}{\familydefault}{\mddefault}{\updefault}{\color[rgb]{0,0,0}$(2,4)$}%
}}}}
\put(1306,-1230){\makebox(0,0)[rb]{\smash{{\SetFigFont{8}{9.6}{\familydefault}{\mddefault}{\updefault}{\color[rgb]{0,0,0}$2$}%
}}}}
\end{picture}%

%% file: pics/CH1.tex
\begin{picture}(0,0)%
\includegraphics{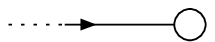}%
\end{picture}%
\setlength{\unitlength}{3947sp}%
\begingroup\makeatletter\ifx\SetFigFont\undefined%
\gdef\SetFigFont#1#2#3#4#5{%
  \reset@font\fontsize{#1}{#2pt}%
  \fontfamily{#3}\fontseries{#4}\fontshape{#5}%
  \selectfont}%
\fi\endgroup%
\begin{picture}(995,356)(3289,-444)
\put(4201,-211){\makebox(0,0)[b]{\smash{{\SetFigFont{10}{12.0}{\familydefault}{\mddefault}{\updefault}{\color[rgb]{0,0,0}$v'$}%
}}}}
\end{picture}%

%% file: pics/CH2.tex
\begin{picture}(0,0)%
\includegraphics{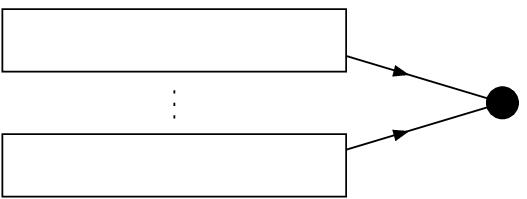}%
\end{picture}%
\setlength{\unitlength}{3947sp}%
\begingroup\makeatletter\ifx\SetFigFont\undefined%
\gdef\SetFigFont#1#2#3#4#5{%
  \reset@font\fontsize{#1}{#2pt}%
  \fontfamily{#3}\fontseries{#4}\fontshape{#5}%
  \selectfont}%
\fi\endgroup%
\begin{picture}(2495,924)(5539,-823)
\put(7951,-211){\makebox(0,0)[b]{\smash{{\SetFigFont{10}{12.0}{\familydefault}{\mddefault}{\updefault}{\color[rgb]{0,0,0}$v'$}%
}}}}
\put(7951,-661){\makebox(0,0)[b]{\smash{{\SetFigFont{10}{12.0}{\familydefault}{\mddefault}{\updefault}{\color[rgb]{0,0,0}$(0,a)$}%
}}}}
\put(6376,-121){\makebox(0,0)[b]{\smash{{\SetFigFont{10}{12.0}{\familydefault}{\mddefault}{\updefault}{\color[rgb]{0,0,0}$N^{\floor}_{d^1,   \bf{k^1}}(\alpha^1  +e_{m^1}, \beta^1)$}%
}}}}
\put(6376,-721){\makebox(0,0)[b]{\smash{{\SetFigFont{10}{12.0}{\familydefault}{\mddefault}{\updefault}{\color[rgb]{0,0,0}$N^{\floor}_{d^t,   \bf{k^t}}(\alpha^t  +e_{m^t}, \beta^t)$}%
}}}}
\end{picture}%

%% file: pics/CH3.tex
\begin{picture}(0,0)%
\includegraphics{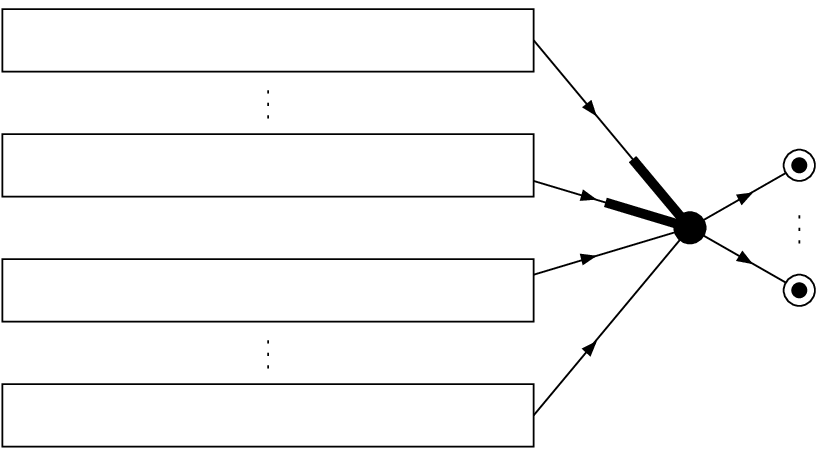}%
\end{picture}%
\setlength{\unitlength}{3947sp}%
\begingroup\makeatletter\ifx\SetFigFont\undefined%
\gdef\SetFigFont#1#2#3#4#5{%
  \reset@font\fontsize{#1}{#2pt}%
  \fontfamily{#3}\fontseries{#4}\fontshape{#5}%
  \selectfont}%
\fi\endgroup%
\begin{picture}(3920,2124)(3514,-3823)
\put(4801,-3121){\makebox(0,0)[b]{\smash{{\SetFigFont{10}{12.0}{\familydefault}{\mddefault}{\updefault}{\color[rgb]{0,0,0}$N^{\floor}_{d^{t'+1},   \bf{k^{t'+1}}}(\alpha^{t'+1} , \beta^{t'+1}+e_{m^{t'+1}}) $}%
}}}}
\put(4801,-2521){\makebox(0,0)[b]{\smash{{\SetFigFont{10}{12.0}{\familydefault}{\mddefault}{\updefault}{\color[rgb]{0,0,0}$N^{\floor}_{d^{t'},   \bf{k^{t'}}}(\alpha^{t'}  +e_{m^{t'}}, \beta^{t'})$}%
}}}}
\put(4801,-1921){\makebox(0,0)[b]{\smash{{\SetFigFont{10}{12.0}{\familydefault}{\mddefault}{\updefault}{\color[rgb]{0,0,0}$N^{\floor}_{d^1,   \bf{k^1}}(\alpha^1  +e_{m^1}, \beta^1)$}%
}}}}
\put(4801,-3721){\makebox(0,0)[b]{\smash{{\SetFigFont{10}{12.0}{\familydefault}{\mddefault}{\updefault}{\color[rgb]{0,0,0}$N^{\floor}_{d^t,   \bf{k^t}}(\alpha^t  , \beta^t+e_{m^t})$}%
}}}}
\put(6846,-2611){\makebox(0,0)[b]{\smash{{\SetFigFont{10}{12.0}{\familydefault}{\mddefault}{\updefault}{\color[rgb]{0,0,0}$v'$}%
}}}}
\put(6846,-3061){\makebox(0,0)[b]{\smash{{\SetFigFont{10}{12.0}{\familydefault}{\mddefault}{\updefault}{\color[rgb]{0,0,0}$(d',a)$}%
}}}}
\end{picture}%

%% file: biblio.tex
\providecommand{\bysame}{\leavevmode\hbox to3em{\hrulefill}\thinspace}
\providecommand{\MR}{\relax\ifhmode\unskip\space\fi MR }

\providecommand{\href}[2]{#2}